`
\documentclass{amsart}

  \usepackage[all]{xy}

  \usepackage{epsf,epsfig,amsfonts,graphicx,color}
  
  \usepackage{url}
   
 \def\R{\mathbb R}

\def\C{\mathbb C}
 \def \E {\mathbb E}
  \def \V {\mathbb V}
   \def \W {\mathbb W}
      \def \S {\mathbb S}

\def\T {\mathbb T}
  
 \def\on{orthonormal }

\def\Ri{Riemannian }

\newtheorem{theorem}{Theorem}[section]
\newtheorem{proposition}{Proposition}[section]

\newtheorem{lemma}{Lemma}[section]
\newtheorem{definition}{Definition}[section]

\newtheorem{example}{Example}[section]

\newtheorem{exer}{Exercise}[section]

\usepackage{amsmath,amssymb}
\usepackage{graphicx}

\begin{document}

\title{Oscillating about coplanarity in the 4 body problem.  }

\author{Richard Montgomery}
\address{Mathematics Department\\ University of California, Santa Cruz\\
Santa Cruz CA 95064}
\email{rmont@ucsc.edu}

\date{October 13, 2018}

\begin{abstract}  For the Newtonian 4-body problem in space we prove that any zero angular momentum
bounded solution suffers infinitely many coplanar instants, that is, times at which  all 4 bodies lie in the same plane. 
This result generalizes a  known  result for collinear instants (``syzygies'')  in the   zero angular momentum  planar 3-body problem,
and extends to  the   $d+1$ body problem  in   $d$-space.  The proof, for $d=3$,  starts  by identifying the
 center-of-mass zero configuration space with real 
 $3 \times 3$ matrices, the coplanar configurations with  matrices whose  determinant  is zero, and the mass metric
 with the Frobenius (standard Euclidean) norm.  Let $S$ denote the signed distance from a matrix to the hypersurface of matrices with determinant zero. 
 The  proof hinges on establishing a harmonic oscillator type ODE
for  $S$ along solutions.  Bounds on  inter-body distances then yield an  explicit   lower  bound $\omega$  for the frequency  of this
oscillator, guaranteeing a degeneration within  every time  interval of length $\pi/\omega$. 
The non-negativity of the  curvature of
  oriented shape space (the quotient of configuration space  by the rotation group)   plays a crucial role in the proof. 
\end{abstract}

\maketitle

\section{Results.} 
Consider the Newtonian  4 body problem in Euclidean 3-space.
Typically,  the four point masses   form the vertices of a tetrahedron.
As the masses move about,  at isolated instants the tetrahedron which they form might  degenerate so that all  
4 bodies  lie on a single  plane.   Must such co-planar instants always occur?

A  solution is called   {\it bounded} if the interparticle distances $r_{ab}$  between the four masses $m_a, a =1,2,3,4$  are bounded for all time
in the solution's domain of definition.  
\begin{theorem}  \label{thmA}  For the 4 body problem in 3-space, 
any bounded zero angular momentum  solution defined on an infinite time interval 
suffers infinitely many coplanar instants. \end{theorem}
  

Theorem \ref{thmA} follows directly from the     finite time interval  oscillation results   of Theorems \ref{thmB} and  \ref{thmC} and Proposition \ref{basic}  below,
results which 
 hold for   the $d+1$-body problem in  $d$-dimensional  Euclidean
 space.   These results generalize the result  \cite{Infinitely}   for the  case $d=2$ of the planar three-body problem .

 Write $q_a \in \R^d,  a = 1, \ldots , d+1$ for
 the positions of the bodies.  Typically,  at each instant  the $q_a$ form the vertices of a    $d+1$-simplex, meaning that 
 their convex hull has nonzero $d$-dimensional volume.  At  special instants 
  this volume  may vanish  by virtue of  all bodies instantaneously  lying  on some affine
 hyperplane.   We call these {\it degeneration instants}.   Write $r_{ab} = |q_a - q_b|$
 for the distances between bodies, $M = \Sigma m_a$ for  the total mass and $G$ 
 for the  universal gravitational constant.  $G$ is included   to get our units straight:
  $GM/{r_{ab}}^3$ has  the  units of $1/(time)^2$, the units of a frequency squared.
 \begin{theorem} \label{thmB} Consider any zero angular momentum solution
 to the standard attracive ($1/r$ potential) Newton's equations for $d+1$ bodies in $d$-dimensional.
 Suppose that along this solution the inter-body distances satisfy the bound 
 \begin{equation}  r_{ab} \le c
 \end{equation}  
  Then, within every   time interval of size $\frac{1}{\pi} (\frac{c^3}{GM})^{1/2}$, 
this solution has a degeneration instant.
 \end{theorem}
{\sc Remark.}   Theorem   \ref{thmB}  represents a quantitative improvement of the syzygy estimates
 found earlier in  the   case $d=2$  described above.  
 
  {\sc Necessity of zero angular momentum in even dimensions.} The regular simplex is a central configuration in all dimensions $d$.  
If the dimension $d$ is even, say $d=2k$,  then one can uniformly rotate    the simplex   in a way consistent with
a splitting of $\R^d$ into $k$ two-planes  to get a relative equilibrium solution to the $d+1$-body problem
in $\R^d$ which has  nonzero  angular momentum
and never degenerates.    These   even-dimensional analogues of the Lagrange rotating equilateral
triangle illustrate  that  for even dimensions $d$  the hypothesis that  the angular momentum be zero
is necessary in theorem \ref{thmB}.  
 
 {\sc General two-body type potentials.}. There is nothing special about the Newtonian $1/r$ potential in  theorem \ref{thmB}.
 It is enough to have a 
  sum of pair  potentials of  the form  
 \begin{equation} V(q) = G \Sigma_{a \ne b} m_a m_b  f_{ab} (r_{ab} (q))
 \label{hyp1}
 \end{equation}
where the  individual two-body potentials $f_{ab}$ are attractive. Specifically, assuming    
 \begin{equation} 
 \label{hyp2}
  f_{ab}' (r) >  0,   f''_{ab} (r)  < 0,  \text{ for } r > 0;  \lim_{r \to \infty } \frac{ f'_{ab} (r)}{r}  = 0.
  \end{equation}
 is enough.  Examples  include the standard Newtonian 3-dimensional gravitational potential   $f_{ab} (r) = - 1/r$ 
 and the power law potentials $f_{ab} (r) = -k_{ab} /r^{\alpha}$ for positive exponent  $\alpha$ and positive constants
 $k_{a b}$.  (We choose the  units  so that  $f_{ab}$ has units $1/(length)$.)
 Hypothesis (\ref{hyp2}) guarantees   that the functions  $f_{ab}' (r)/r$ are positive and   strictly monotone
 decreasing   so that  for each $c>0$ and pair $ab$ we have that 
  $ r_{ab} \le c \implies \frac{f'_{ab} (r_{ab})} {r_{ab}}  \ge \delta_{ab} :=   \frac{f'_{ab} (c)}{c}$. 
 Taking   $ \delta$ to be the minimum of these  $\delta_{ab}$  over all pairs
   we get  
 \begin{equation}
 \label{hyp3}
 r_{ab} \le c  \text{ for all pairs } ab  \implies \frac{1}{r_{ab}} f'_{ab} (r_{ab}) \ge \delta >  0  \text{ for all pairs } ab.
  \end{equation}   
  Then, we have
   \begin{theorem}  \label{thmC} Consider the  zero angular momentum   Newton's equations for $N= d+1$ bodies 
   moving in Euclidean  $d$-dimensional space 
    under the influence of  the attractive potential
 (\ref{hyp1}) whose  2-body potentials satisfy hypothesis (\ref{hyp2}).    Suppose
 that  along such a   solution  all 
 its inter-body  distances $r_{ab}$  satisfy the bound $r_{ab} \le c$. Then,  in every   time interval of size $(GM \delta)^{-1/2}/\pi$, 
this solution has a degeneration instant.  Here  
 $\delta$ is as in implication (\ref{hyp3})  above,  and $M$ the total mass.
 \end{theorem}
 
 We now describe the  key ingredients behind these   Theorems. 
\begin{definition} $\Sigma$ is  the degeneration locus within configuration space -- the set of configurations for which the $d+1$ masses all lie on
 a single affine hyperplane. 
 \end{definition}  $\Sigma$ is a (singular) hypersurface in the full configuration space which 
 cuts it into two disjoint congruent halves, the   simplices having  positive  volume, and those having  negative volume.  (The sign  
 of the volume depends on the  orientation of Euclidean space and the ordering of the masses, which we fix once and for all. ) 
 Write $sgn(det(q))$ for the sign of the volume,   defined for $q \notin \Sigma$.  
 For example, if  $d=3$,  then  $sgn(det(q))$ is the sign of the triple product
 $(q_2 -q_1) \cdot ((q_3 -q_1) \times (q_4 -q_1))$.  
  
 \begin{definition} The signed distance $S(q)$ of a configuration $q$ of $d+1$ point masses in $\R^d$
is the  distance from $q$ to the degeneration locus relative to the mass inner product (described in  subsection \ref{ssNewton}),
 that distance being given a plus sign if   the signed volume of $q$ is positive
and a minus sign if negative. In symbols:
$$S(q) = sign(det(q))dist(q, \Sigma).$$
with $S(q) = 0$ if and only if $q \in \Sigma$.
\end{definition}
 In Prop. \ref{propSVD} below   we prove that $|S(q)|$ is the smallest singular value of a $d \times d$ matrix 
 representing $q$ in the center-of-mass frame.

\begin{proposition}
\label{basic}  [Main computation.]  If $S$ is smooth along a  zero angular momentum solution 
$q(t)$  to 
Newton's equations  then $S(t): = S(q(t))$ evolves according
to $$\ddot S = - S g (q, \dot q), \text{ with } g >0 \text{ everywhere } .$$
If, moreover,   all interparticle distances $r_{ab}$ satisfy   $r_{ab} \le c$
then  $g \ge GM/c^3$ for  the Newtonian ($f_{ab} (r) = -1/r$) potential  case,
and,  more generally, $g \ge GM \delta$  for potentials of the form (\ref{hyp1}), with (\ref{hyp2})  in force  and  $\delta$ as per (\ref{hyp3}).
\end{proposition}

 \textbf{Acknowledgements: } I would like to thank Alain Albouy,  Gil Bor,   Joseph Gerver, Connor Jackman, 
Adrian Mauricio Escobar Ruiz,   Robert  Littlejohn, and Rick Moeckel  for useful discussions. This material is based upon work supported by the National Science Foundation under Grant No. DMS-1440140 while the author was in residence at the Mathematical Sciences Research Institute in Berkeley, California, during the Fall 2018 semester.
  
 \section{Motivation and Main ideas.}
 
 Newton's  N-body equations in d-space are invariant under the  isometry group of the inertial Euclidean space, $\R^d$,
 so we can push them down to form  a system of ODEs on   ``shape space'', by which we mean the quotient space
 of the N-body configuration space by the isometry group  of $\R^d$.  There are actually two shape spaces, depending on whether or not 
 we allow orientation
 reversing isometries.  In the body of this paper we will    work  on  the ``oriented shapes space'' which arises from  taking the quotient with respect the group $SE(d)$ of  orientation-preserving isometries.  In   appendix B
 we describe the relationships between these two  shape spaces.
   
  We will speak of ``downstairs'' to mean  we are working  on  the quotient 
  and ``upstairs''
  to mean  we are working on the original  configuration space.  Upstairs, Newton's equations have the form 
$\ddot q = -\nabla V (q)$.  Downstairs, the equations   have precisely  this same form 
  {\bf provided}  that  the
 total angular momentum is zero \footnote{If the angular momentum is non-zero there  are additional `magnetic'
 terms'' in the equations downstairs, meaning terms   linear in   velocities,  and also additional equations
 involving   `internal
 variables'' which represent   instantaneous rigid body tumbling coupled to dynamics on the shape space,
 these internal variables lying in   co-adjoint orbits for   $SO(d)$.}.
 In  writing down the downstairs  zero-angular momentum Newton's equations,   the  acceleration
 $\ddot q$ is   replaced by the   covariant acceleration $\nabla_{\dot q} \dot q$  where $\nabla$ is the
 Levi-Civita connection arising out of the   the induced {\it shape metric} downstairs.
 This  shape metric, induced by the flat   kinetic energy metric upstairs, is curved.  

 Robert Littlejohn  \cite{Littlejohn} pointed out to me  that  the oriented shape space for
 $d+1$ bodies in $\R^d$ is  homeomorphic to a Euclidean space. 
This  topological fact is  well known  for  the case  $d=2$ of the 3 body problem  in the plane where it  has 
 proven to be of great utility.   For
higher $d$ the  fact has been known for some time amongst certain statisticians
and   can be found  in  \cite{Kendall} and \cite{Le} .  Although I do not use this
fact here, it is this single 
fact that inspired my faith that something like the theorems in this paper might hold.

Not only is the oriented shape space homeomorphic to a Euclidean space, but it
is smooth at most points.    The points where it fails to be smooth are those shapes of 
corank $2$ or higher. (The locus of such points has codimension $4$.) Here we use the following terminology 
\begin{definition} The corank of a configuration $q = (q_1, \ldots , q_{d+1})$, or of its corresponding  shape, is the codimension
of the smallest affine subspace in  $\R^d$  which contains   all $d+1$  of the vertices $q_1,  \ldots , q_{d+1}$.   The {\it rank} of 
a configuration is the dimension of this smallest affine subspace.  
\end{definition}

We continue to write $\Sigma$ for the degeneration locus, either upstairs or downstairs.
Downstairs, in {\it oriented} shape space,  $\Sigma$ is   a   totally geodesic hypersurface, at least at its  smooth points.   (This is a bit strange 
since $\Sigma$  is not totally geodesic upstairs. For example, when  $d=3$, imagine  connecting two quadrilaterals which lie in different planes
by geodesics, i.e.  straight lines between vertices.  The resulting curve in configuration space  is non-degenerate at most instants.)  
 To see the total geodesy downstairs, select any orientation reversing isometry $R$, for example, in the case $d=3$,
a reflection about the xy plane.    Downstairs  $R$ is an isometry of shape space 
whose fixed point set is precisely $\Sigma$.  A general theorem in Riemannian geometry now  implies that
$\Sigma$ is   totally geodesic.  In terms of Newton's equations, this `total geodesy' is basically the assertion that, for example for the case  $d=3$,  that a configuration
which initially lies in a plane, and all of whose velocities initially  lie in that plane, will remain in that plane for all time.

{\sc Heuristics.} Proposition \ref{basic}  
asserts that the signed distance $S$  from $\Sigma$ behaves qualitatively  like a one-dimensional  harmonic oscillator, oscillating  around $S = 0$. 
The  physical intuition behind this  phenomenon was pointed out to me by Mark Levi many years ago.
The potential is invariant under isometries so descends to a function downstairs. How to interpret  this potential
downstairs?  Write $\Sigma_{ab} \subset \Sigma$ for the binary collision locus
$r_{ab} = 0$.  One computes that $r_{ab} (s)  = \mu_{ab} dist(s, \Sigma_{ab})$ where $dist(s, \Sigma_{ab})$ is the distance
from $s$ to $\Sigma_{ab}$ and 
where $\mu_{ab} = \sqrt{M / m_a  m_b}$.  
Consequently, re-interpreted  downstairs,    formula (\ref{hyp1})  for the
potential asserts that a  point $s$ in shape space is subjected to the force
of an  attractive potential exerted by the $d \choose 2$  sources  $\Sigma_{ab} $, all of which lie in the  
``hyperplane'' $\Sigma$.  So, of course, the shape  is always attracted to $\Sigma$! 
And  as long as the shape's  ``vertical' kinetic energy is not too large, it will always return to cross 
$\Sigma$,  
oscillating  forever back and forth across the attracting  `hyperplane' $\Sigma$.  

{\sc Choice of  S  versus signed volume. }
In \cite{Infinitely}, in proving  Theorem \ref{thmA} for the case $d=2$, I used a function $z$ in place of the  $S$ of proposition \ref{basic}.
This  $z$ was  the signed area of the oriented  triangle normalized by 
divided it by the   moment of inertia $I$ that the triangle  would have   if all masses were assigned the value  $1$.     The obvious generalizations $z_d$  of  this  $z$ to $d > 2$, namely a normalized signed volume, did  not work out.  All my attempts at proving a version of   proposition \ref{basic}  for such a function in place of $S$ failed.
The function $z_2$  satisfies a kind of  monotonicity relation with respect to  geodesics  orthogonal to $\Sigma$ which fails for $z_d, d > 2$ and this monotonicity
was required to get positivity  of $g$ in  Proposition \ref{basic}.  
The need for such a relation  led to introducing $S$.  After the fact, one  observes that
the identity  $z = S/\sqrt{I}$
holds   for equal masses when $d=2$, and  fails for $d > 2$.)

 {\sc Key ingredients to the proof.}  The proof of Proposition \ref{basic} relies on four 
 key facts .  
 \begin{itemize} \item{}Fact 1.  $S$ satisfies the Hamilton-Jacobi equation $\| \nabla S \| = 1$
 wherever $S$ is smooth. This fact implies that   the integral curves of the  gradient flow of $S$ are geodesics.
 \item{}Fact 2.  The shape metric   is everywhere non-negatively curved.  
 \item{}Fact 3.   There is a close  relationship between the sign of the  second fundamental form of 
   distance level sets  ( the $\{ S = t \}$'s)  from a totally geodesic submanifold ($\Sigma= \{S = 0 \}$ ) and   the sign of the    curvature of the ambient space within which
the level sets lie.    
 This relation is detailed on p. 34-37 of Gromov \cite{Gromov} and recalled below as proposition \ref{curv}. 
 
 \item{} Fact 4.  (Theorem \ref{SSSVD}).  The singular locus of $S$ has codimension $2$.  This
 locus, denoted $Sing(S)$ below, consists of all points at which $S$ is not smooth.
 \end{itemize}

\section{Set-up and Reduction.} 

The proofs of all of our  theorems hinge on Proposition  \ref{basic} which is  a computation.
We achieve the computation  by exploiting the relations  between Newton's equations at zero angular momentum as expressed  upstairs on the usual
configuration space and downstairs on shape space.  The process of pushing the equations downstairs is 
referred to as ``reduction''. Our reduction procedure is a  metric reduction, putting kinetic energy to the fore, 
as opposed to the oft-used symplectic reduction.  The two reduction procedures are formally equivalent but the metric approach makes
out computation tractable.    In this 
section we go through the reduction  for  the case $d=3$.  At the end,   in subsection \ref{d_general}, 
we   describe the   small changes needed for the set-up of reduction for higher $d$.  

Write $M(k,m)$ for the space of $k \times m$ real matrices.  The configuration space for the
  4 body problem in $\R^3$ can be naturally identified with the space 
  $M(3,4)$. To do so,   think of the  four  vectors $q_1, q_2,q_3, q_4 \in \R^3$ which
   define the positions of the four bodies  as   column vectors and place them side-by-side to form the    $3 \times 4$ matrix 
\begin{equation}
q = \left(
\begin{array}{cccc}
 q_1 & q_2  &  q_3 & q_4 
\end{array}
\right)  \in M(3, 4).
\label{config_rep}
\end{equation}


The translation subgroup $\R^3$  acts on $M(3,4)$ by  $q_a \mapsto q_a +b$,  $b \in \R^3$, $a = 1,2,3,4$, which in matrix terms is
\begin{equation}
q \mapsto q +  (b,b,b,b)
\label{transl}
\end{equation}
 The quotient  of  $M(3,4)$ by this action
can be identified with the matrix space  $M(3,3)$.  This identification depends on choosing a basis for
the 3-dimensional subspace $x_1 + x_2 + x_3 + x_4 = 0$ of the mass label space $\R^4$, or,
what is the same thing, a choice of ``Jacobi vectors''. The results are independent of this choice
\footnote{Usually the   quotient $M(3,4)/\R^3$  is identified with   the codimension 3 linear subspace
of $M(3,4)$ obtained by   fixing the center of mass to be zero. There is a {\it mass independent} way to form the  identification of the quotient  with $M(3,3)$ which 
is more useful for our purposes. This alternative  perspective, due to Albouy and Chenciner REF, is  is reviewed in an appendix.
It is equivalent to fixing the center-of-mass, once
masses  are chosen.}.
See Appendix  A for details.  The rank of a configuration becomes the rank of the representing  matrix  so that  the degeneration locus is 
$$\Sigma = \{q \in M(3,3): det(q) = 0 \}.$$

\subsection{Oriented Shape Space.}
\label{ssShapeSpace}
Rotations  of $\R^3$   act on both $M(3,4)$ and its translation-quotient $M(3,3)$ by 
$$q \mapsto g q ,  g \in SO(3). $$
\begin{definition} 
The oriented shape space $Sh= Sh(3,4)$ is the  topological quotient space  $M(3,3)/SO(3)$.  The  
 quotient map $M(3,3) \to Sh(3,4)$ will be denoted by $\pi$.   The projection of a configuration $q \in M(3,3)$
 will  be called ``the shape'' of $q$.  
\end{definition} 

{\sc Remark.}  The group of orientation-preserving isometries of $\R^3$, denoted $SE(3)$,  is made up of the  translations ($\R^3$)  and 
the rotations ($SO(3)$). 
We can naturally identity $Sh$ with the quotient $M(3,4)/SE(3)$, by using reduction in stages:  first quotient by translations $\R^3$
to get $M(3,3)$ and then by rotations $SO(3)$ to get to $Sh$.    The projection $M(3,4) \to Sh$
will also be denoted by $\pi$.  

Recall that the action of a group $G$ on a set $Q$ is called free if $gq = q \implies g = id$.
It is well-known (see for example Prop 4.1.23 of \cite{Abraham}) that the free action of a compact Lie group $G$ on a smooth manifold $Q$
yields a quotient $Q/G$ which is itself a smooth manifold with the quotient map $Q \to Q/G$ being a smooth
projection. 
The action of $SO(3)$ on $M(3,3)$ is free on the open dense subset of $M(3,3)$
consisting of matrices of rank 2 and 3, i.e. on the planar and spatial configurations.  
Moreover, $M(3,3)$ is stratified by rank. The rank $3$ matrices form an open dense subset whose complement is the singular hypersurface $\Sigma$.
The  rank $2$ matrices form the generic points of $\Sigma$ , the points at which it is smooth.  The rank 1 points, i.e. the collinear configuratins,
have   codimension $4 = 2*2$ in $M(3,3)$. See for example,  \cite{Arnold} for this computation. There is only one   rank $0$ matrices, namely the $0$ matrix which represents total collision.   . 
 Hence we get 

\begin{proposition}
\label{shape1} Let $U \subset M(3,3)$ denote the set  of rank 2 and 3 matrices, henceforth referred to as ``generic configurations''. 
  Then the restriction of 
 $\pi: M(3,3) \to Sh(3,4)$  to $U$ gives the space of  rank 2 and 3 shapes within $Sh(3,4)$  a smooth structure
 in such a way that this restricted projection is   a smooth submersion. Moreover this restricted projection  $\pi: U \to \pi(U)$  gives $U$ the structure of a principal $SO(3)$ bundle.   The complement of $U$ has   codimension 4 within  $M(3,3)$. 
\end{proposition}

\subsubsection{Newton's Equations.} 
\label{ssNewton}To write down Newton's equations for the motion of the 4 bodies,  we need the potential and the choice of masses.
We have   written down the  potential (eq (\ref{hyp1})). 
The choice of  masses $m_a > 0$ for each body defines an inner product $\langle \cdot , \cdot \rangle$ on $M(3,4)$ called the ``mass metric'' or ``kinetic energy metric''
according to  
$$\frac{1}{2} \langle \dot q , \dot q \rangle  = \frac{1}{2} \Sigma_{a=1} ^4 m_a | \dot q_a |^2 .$$
We use    the absolute
value symbol  for the usual norm in our Euclidean inertial $\R^3$.
When we interpret $\dot q \in M(3,4)$
to represent the velocities of the four bodies then   the above expression is the
usual expression for the total kinetic energy.     Newton's equations can now be written 
\begin{equation}
\label{Newton}
\ddot q = -\nabla V (q)
\end{equation}
 where the gradient $\nabla V$ is computed using the mass inner product:
 $dV(q) (\delta q) = \langle \nabla V(q), \delta q \rangle$.  
 
 The mass inner product induces an inner product  on the
  translation reduced configuration space $M(3,3)$ by declaring the projection $M(3,4) \to M(3,3)$
  to be a metric projection.  Equivalently, we can view $M(3,3)$ as
  a subspace of $M(3,4)$ by fixing the center of mass to be zero, and then take the restricted inner product. (Again, see   appendix A.) 
  We can choose a basis for the   1st $\R^3$ factor of
  $M(3,3) = Hom(\R^3, \R^3)$ such that the inner product becomes
  $$\langle q, q \rangle = Tr(q^t q ), $$
 namely, the inner product on matrices  for which the matrix entries form an orthonormal
  linear coordinate system.  Such a basis  corresponds  to a choice
  of {\it normalized Jacobi vectors}.  See Appendix A.  
  So done, Newton's equations have precisely the   form, eq (\ref{Newton}) when written on $M(3,3)$.
  
  \subsection{Reduced Newton's equations.}
  \label{ssReducedNewton}
  We  push Newton's equations  and the kinetic energy metric down to shape space.
  For this purpose it will be helpful to keep in mind the following generalities. 
  
  {\sc Metric projections and Riemannian submersions.}  Whenever we have a metric
  space $M$ with distance function $d_M$ and an onto map $\pi: M \to B$ we can {\it try} to define a metric $d_B$  on $B$
  by $d_B (b_1, b_2) = d_M (\pi^{-1} (b_1), \pi^{-1} (b_2))$, or, in English, the distance between points downstairs is the distance
  between their corresponding fibers upstairs.  When this construction
  works we say that  $\pi:M \to B$ is a metric projection or {\it submetry}.   If  $B = M/G$ is the quotient of $M$
  by the action of a compact Lie group acting on $M$ by isometries  and $\pi$ is the quotient projection then the construction always
  works.  
  If, in addition, $M$ is a manifold whose  metric $d_M$ comes from a \Ri metric and if  the 
$G$-action  is free so that the quotient map $\pi$ is a smooth submersion with  smooth  $B$,  then the induced   distance function  $d_B$ also
arises as   the distance function of a   \Ri metric on $M$.   In this case  $\pi:M \to B$ is a  {\it \Ri submersion} which has the following
infinitesimal meaning.   The `vertical space''
  $V_q \subset T_q M$ through $q \in M$  is  defined to be the kernel of $d \pi_q$; equivalently, it is the tangent space at $q$  to the fiber $\pi^{-1} (s) = Gq$ through $q$. Define the ``horizontal space'' $H_q$ to be  the orthogonal complement to the vertical:   $H_q = V_q ^{\perp}$.  Then  the restriction of $ d \pi_q$ to $H_q$ is a linear isomorphism.  
  Declaring this   linear isomorphism $H_q \to T_s B$ to be an isometry  induces an  inner product on $T_s B$, and this inner product  is 
  independent of the point $q \in \pi^{-1} (s)$ since $G$ acts isometrically on $M$.    
 Distance minimizers between fibers upstairs are geodesics in $M$ orthogonal to the fibers.  From this
  follows the well-known fact that geodesics orthogonal to fibers at one point are orthogonal at every point,
  and that the geodesics downstairs in $B$ are precisely the projections of horizontal geodesics upstairs.
  
  \vskip .4cm
  
  In this way, starting from the mass metric on $M(3,4)$ or $M(3,3)$,  we get a   metric on $Sh= Sh(3,4)$ which is \Ri at the generic shapes
  (those of rank 2 or 3)  and over these points   is such that $\pi: M(3,3) \to Sh$ is a \Ri submersion.    
  
    \vskip .4cm
    
    To push Newton's equations down to $Sh$ we must understand the {\it dynamical meaning}
    of being horizontal in $M(3,3)$.  In \cite{isoholonomic}  (or  \cite{tour}) I compute  that $\dot q$ is orthogonal to the $SO(3)$
    orbit through $q$ if and only if the total angular momentum  $J(q, \dot q)$ of the pair $(q, \dot q)$ is zero.
    The expression for $J$ as a function on $TM(3,4) = M(3,4) \times M(3,4)$ is 
   $J ( q, \dot q ) = \Sigma m_a q_a \wedge \dot q_a$
   for $M(3,4)$, and   is the same when restricted to  $TM(3,3)$ viewed as  subspace of $TM(3,4)$.    Recall that $J$ is conserved for any
  potential of the form of eq (\ref{hyp1}), that is to say  $J(q(t) , \dot q (t)) = J(q(0), \dot q (0))$ along solutions $q(t)$ to Newton's equations.
  Now let  $\nabla$ be the Levi-Civita connection for the shape metric. Observe that since the potential is $SE(3)$
  invariant it also defines a projection on $Sh$.  We will use the same symbol $V$ for the potential upstairs and downstairs. 
  We have 
  \begin{lemma} 
  \label{reducedNewton}
 Any zero angular momentum solution  to    
 Newton's equations passing through generic (i.e. rank 2 and 3)   points of $M(3,3)$ projects to a  curve $\gamma$ in shape space which satisfies  
 $$\nabla _{\dot \gamma} \dot \gamma = - \nabla V (\gamma(t)). $$
Conversely, the horizontal lift  of any such solution is a   zero-angular momentum solution to Newton's equations  upstairs.
\end{lemma} 
Regarding   `horizontal lift'' see, again \cite{isoholonomic} or  chapter 13 of \cite{tour}.

{\sc Proof.} This theorem is a general fact, holding for any Hamiltonian of the form kinetic plus potential on any manifold endowed
with the smooth free action of a Lie group which keeps both the kinetic (metric) and potentials invariant.
For a proof  see for example, \cite{tour}.  
 
 The special case when $V =0$ will be useful below.
 \begin{lemma} 
 \label{lines} 
 Any zero angular momentum straight line $q + tv$ in $M(3,3)$   projects to a geodesic in Shape space $Sh(3,4)$.  
Conversely, the horizontal lift  of any geodesic in $Sh(3,4)$   is a   zero-angular momentum straight line in $M(3,3)$.   The geodesic is parameterized
by arc length if and only if $\|v \| = 1$. 
\end{lemma} 

\subsection{Set-up for general dimension $d$}\label{generalD}
\label{d_general}

{\sc Going from $d=3$ to general $d$.}   The configuration space for $N$ bodies in $\R^d$ is the space $M(d, N)$ of $d \times N$
real matrices.  Its  quotient
by the translation group  of $\R^d$ forms an  $M(d, N-1)$ once a basis for the hypersurface $x_1 + x_2 + \ldots + x_N = 0$
of the  mass-label is chosen.  Again see   Appendix A.  In our case of $N = d+1$ we thus get the translation-reduced
configuration space  $M(d,d)$   of square matrices.   The degeneration locus $\Sigma$  is  given by  $\{q: det(q) = 0\}$.  Shape space is $Sh(d,d+1) = M(d,d)/SO(d) = M(d, d+1)/SO(d)$.  Proposition
\ref{shape1} holds with `rank $2$ and $3$' replaced by `rank $d-1$ and rank $d$.

Introducing masses puts an inner product on the mass label space, and so  on $M(d,d+1)$ and on its translation quotient   $M(d,d)$.
The masses also allow us to identify $M(d,d)$ as a linear  subspace, rather than a quotient space,  of $M(d, d+1)$, 
namely as the subspace of  center-of-mass zero configurations.   An \on basis for the  hypersurface $\Sigma x_i = 0$ 
is equivalent to a choice of  normalized Jacobi vectors and relative to these coordinates the   mass-induced  
 inner product structure on $M(d,d)$ is      standard :   $\langle q, q \rangle  = tr (q^t q) = \Sigma_{i,j} q_{ij}^2$,
 (This is the inner product whose associated norm is called the Frobenius norm).
 Relative to this inner product the  rotation group $SO(d)$ acts isometrically by left multiplication and the action leaves  the potential invariant and so the metric and the zero-angular momentum Newton's 
 equation push down to the quotient shape  space $Sh(d, d+1) = M(d,d)/SO(d)$. 
 The reduction lemmas \ref{reducedNewton} and \ref{lines}  for this reduced dynamics hold as stated, upon  replacing `3' by `$d$' in the obvious places.

\section{ The prooof of  prop \ref{basic}: signed distance as an oscillator.} 

We proceed to differentiate $S$ along a solution arc which does not pass
through  any singular point  of $S$. We have 
$$\dot S = \langle \nabla S, \dot \gamma \rangle$$
so that
\begin{eqnarray}
\label{Sdot} 
\ddot S & = \langle \nabla S, \nabla_{\dot \gamma} \dot \gamma  \rangle
+  \langle  \nabla_{\dot \gamma} \nabla S,  \dot \gamma  \rangle \\
& = \langle \nabla S, - \nabla V  \rangle
+  \langle  \nabla_{\dot \gamma} \nabla S,  \dot \gamma  \rangle
\end{eqnarray}
We    estimate each term of this last equation separately, showing that each term has  the
form $-S g$ with $g \ge 0$.  We  verify that the `$g$'  for the first term  is always positive and satisfies the stated bounds when $r_{ab} \le c$.  

{\bf  First term, } $\langle \nabla S, - \nabla V  \rangle$. At smooth points of $S$, the integral curves of the vector field $\nabla S$
are geodesics orthogonal to the level sets of $S$, and in particular to the level set $S=0$ which is  the degeneration locus $\Sigma$. 
This fact holds true generally for  the   signed distance function from a hypersurface on
any Riemannian manifold, and  is      closely related
to the fact that   signed distance  satisfies the Hamilton -Jacobi equation:  
$\| \nabla S \| =1$. 
 
 We proceed in the special case of $d=3$  for this paragraph,  for simplicity.  The   geodesics in $M(3,3)$, or
 in shape space,  are the projections
  of straight lines $q + t v$  in $M(3,4)$ for which  $(q,v) \in M(3,4) \times M(3,4)$ has zero total angular
momentum and zero total linear momentum. See lemma \ref{lines} above. 
(Zero linear momentum arises from working on $M(3,3) \subset M(3,4)$
by identifying it with the zero center-of-mass configurations. Alternatively, having zero linear momentum is equivalent to the assertion that the velocity $v$
is orthogonal to the translation action.) 
 The parameter $t$ is arclength
provided $\langle v, v \rangle = 1$. Now the smooth points $q$  of the degeneration
locus $\Sigma$ are the planar points.  In order for a  geodesic
to be perpindicular to $\Sigma$ at such a  $q$ we must have
that  $v$ is    perpindicular to all $\delta q \in  T_q \Sigma$.
By rotating, we  may   assume that the 4 vertices of $q$ lie in the $xy$ plane which we will denote by ``$\R^2$''. 
Then any  variation $\delta q = (\delta q_1, \delta q_2, \delta q_3, \delta q_4)$
with $\delta q_a \in \R^2$   represents a planar variation of $q$ and hence a  tangent vector to $\Sigma$ at $q$. Since
$$\langle \delta q, v \rangle = \Sigma m_a (\delta q_a) \cdot v_a$$
and since the $\delta q_a$ are arbitrary vectors in $\R^2$, we see that our tangent vector $v$ must have all 4 of its component
vectors $v_a$ perpindicular to $\R^2$, which is to say, along the z-axis.
But then, along our geodesic the squared inter-body distances are  \begin{equation}
\label{squaredDistances}
r_{ab} ^2 = | q_a + t v_a ) - (q_b + t v_b ) |^2 = r_{ab} (0)^2 + t^2 | v_a - v_b |^2
\end{equation}
where the cross term is zero since $q_a, q_b$ lie in $\R^2$ while $v_a, v_b$
are orthogonal to $\R^2$.  

For general $d$, equation (\ref{squaredDistances}) continues to hold for
a geodesic orthogonal to the degeneration locus.  Indeed, the only real   difference
between  the proof above for $d=3$ and the proof for $d>3$ is notational.
Now the   $q_a$,  representing a point on the degeneration locus,   can be taken to all  lie in  a fixed affine hyperplane 
of $\R^d$ so that the variations   $\delta q_a$, $a =1, \ldots  N = d+1$ 
can be taken to be arbitrary vectors tangent to the correspoding linear hyperplane $\R^{d-1}$.  As a consequence the  
$v_a$   all lie in the one-dimensional orthogonal to this $\R^{d-1}$ and the computation is the same.

Now look at the negative of the  potential in the gravitational case:
$$U = - V = G \Sigma \frac{m_a m_b} {r_{ab} }$$
along our geodesic. Each individual term  $\frac{m_a m_b} {r_{ab}}$ is 
{\it strictly decreasing} or constant in  $t^2$.   
Indeed $\frac{d}{dt} \frac{1}{r_{ab} (t)} = - \frac{t |v_{ab}|^2} {r_{ab} (t) ^3} = - S \frac{ |v_{ab}|^2} {r_{ab} (t) ^3}$, since $S = t$ as long as the geodesic is the  unique minimizer to the degeneration locus. Summing, we obtain
$$ \langle \nabla S, - \nabla V  \rangle = \langle \nabla S,  \nabla U  \rangle = - S g_1$$
with
$$g_1 =  G \Sigma m_a m_b \frac{ |v_{ab}|^2} {r_{ab} (t) ^3} > 0$$ 
as desired. 

If each $r_{ab}$ is bounded above by $c$, we have that  $g_1 \ge \frac{G}{c^3} \Sigma m_a m_b |v_{ab}|^2 $.
But, if $\Sigma m_a v_a = 0$, 
we find  that $\| v \|^2 =  \Sigma m_a m_b  |v_{ab}|^2 / M$ (``Lagrange's identity'') 
and  since we have that $\|v \|^2 = 1$  (since $t$ is arclength)   it follows that $\Sigma m_a m_b |v_{ab}|^2 = M$  which yields
$g_1 \ge GM /c^3$, which completes the proof for the gravitational case.  

In the case of a general potential satisfying hypothesis (\ref{hyp1}), (\ref{hyp2})  we get that
$\frac{d}{dt} f_{ab}(r_{ab}) = f_{ab}'(r_{ab}) \frac{d}{dt} (r_{ab} (t)) = f_{ab}'(r_{ab}) (t v_{ab}^2)/r_{ab} = S (\frac{f'(r_{ab})}{r_{ab}}) (v_{ab})^2$.
Summing, we get $\langle \nabla S, - \nabla V  \rangle =  - S g_1$ with 
$g_1 = G \Sigma m_a m_b (\frac{f'(r_{ab})}{r_{ab}}) (v_{ab})^2 >0 $.  Under the boundedness
assumption, eq (\ref{hyp3}) yields that   $\frac{f'(r_{ab})}{r_{ab}} > \delta$ for all pairs $a,b$
and the lower bound for $g_1$   proceeds exactly as in   the previous paragraph.  

QED for Term 1. 

{\bf Second term,}$\langle \nabla_ v \nabla S,v \rangle$. 
 For  a  fixed shape  $p$,  $p \notin Sing(S)$   
  $$v \mapsto Q_p (v,v) := \langle \nabla_ v \nabla S,v \rangle, v \in T_p Sh$$  
 is a quadratic form on the tangent space $T_p Sh$.  
{\bf We will show that }$Q_p (v,v) = -S(p) H_p (v,v)$ {\bf  where} $H_p \ge 0$ {\bf is a positive semi-definite quadratic form. }
The trick for achieving this inequality  is to recognize the quadratic form $Q_p$ 
as being essentially  the second fundamental form   of the equidistant hypersurface $\Sigma_t$
 from $\Sigma$ which passes  through $p$, namely  
$$\Sigma_t : = \{S = t\}; \text{ where }  t= S(p)$$ and then 
to use a  relation 
 between  the   sign of such 
second fundamental forms and the sign of   the   ambient   curvature.  

Take $v = \nabla S$ in $Q_p (v,v)$.   Differentiate the identity  $\langle \nabla S,  \nabla S  \rangle =1$
with respect to $v$ 
to see that $ \langle \nabla_v \nabla S, v \rangle = 0$, so that  $Q_p (v, v) = 0$.

Take $v \perp \nabla S$. Then $v$ is 
tangent to $\Sigma_t$ while $\nabla S $ is the unit normal $N$  to $\Sigma_t$. 
Recall that  second fundamental form to a hypersurface $V$  with  unit normal vector field $N$ is the quadratic form
$\Pi(v,v)= v \mapsto \langle \nabla_v N , v \rangle$ defined for vectors  $v$ tangent to $V$.  It follows that
for   $Q_p (v, v) = \Pi_p (v,v)$ for $v \perp \nabla S$ is the   second fundamental form $\Pi_p$ of  the hypersurface $\Sigma_t$ at the point $p \in \Sigma_t$.
Summarizing:
$$Q_p (v,v) = \begin{cases} 0 \text{ for }  v \parallel \nabla S \\
                     \Pi_p (v,v) \text{ for }   v \perp \nabla S 
                     \end{cases}
                     $$

 We recal  some facts about the second fundamental form $\Pi$   of a hypersurface.  
\begin{itemize}
\item{}(1) A hypersurface is {\it totally geodesic} if and only if $\Pi = 0$. 
\item{}(2)Replacing the choice of  unit normal $N$ to the hypersurface by
its negative $-N$ replaces  $\Pi$ by its negative $- \Pi$.  
\end{itemize}

Our hypersurface $\Sigma$ is totally geodesic, as mentioned earlier in `heuristics'.
Indeed,  $\Sigma$ is the fixed point set of an isometric involution $i: Sh \to Sh$
and  fixed point sets of isometric involutions are always totally geodesic.  This isometric
involution $i$, called  ``reflection about $\Sigma$'',    is
implemented by the nontrivial element of the  two-element group $O(d)/SO(d)$. 
Any orientation reversing orthogonal transformation $R \in O(d)$
realizes this nontrivial element and acts on   shape space 
by sending the shape $s = \pi(q)$ to $i(s) = \pi(Rq)$. 
Now $i^* S = -S$ from which it follows that $i_* \nabla S = - \nabla S$, and thus, using item (2) above, that 
$i^*Q = - Q$.  It follows that we can write  $Q = -S H$ where
$i^*H = H$.  
 It remains to show that $H$ is positive semi-definite.
 
 \vskip .3cm

{\sc A necessary detour into curvatures.}

\begin{definition} A hypersurface  is {\it convex} relative  to the choice of normal  $N$ if $\Pi \ge 0$ for this choice of normal,
and {\it concave} relative to $N$  if $\Pi \le 0$ for this choice of normal. 
\end{definition}

\begin{example} The boundary of a  convex domain having  smooth boundary in Euclidean space is convex in the above sense provided we
use the outward pointing normal.  
\end{example}

Let $M$ be a \Ri manifold and $V \subset M$ a   hypersurface in $M$,   together with  a choice of unit normal  $N$ along the hypersurface.
Then close to $V$ we  have the family $V_s,  -\epsilon < s < \epsilon$ 
of nearby equidistant hypersurfaces   formed by travelling along  the   geodesics tangent to the unit normal  $N$
for a distance $s$.   By flowing along these geodesics we also have  diffeomorphisms $$\phi_s:  V \to V_s.$$
 Write $\Pi_0$ for the second fundamental form of $V$ relative to $N$ and $\Pi_s$ for that of the equidistant $V_s$. 
 Recall that we say  that $M$ is ``non-negatively curved'' if its sectional curvatures are all positive or zero,
 and ``non-positively curved'' if all of its sectional curvatures are all negative or zero. 
The following  basic  relationship between extrinsic and intrinsic curvature.
  is  found on p. 34-37 of \cite{Gromov}. 
\begin{proposition} 
\label{curv} (See figures 1).  If  the ambient curvature of the Riemannian manifold $M$  is non-negative
and if the hypersurface  $V \subset M$ is   concave with respect to  the choice of unit normal $N$ for $V$, then its positive equidistants
$V_s, s> 0$ are at least as concave as $V$:
$\phi_s ^* \Pi_s \le \Pi_0 \le 0$ for $s > 0$. 

If the ambient curvature of $M$  is non-positive 
and if the hypersurface  $V \subset M$ is   convex  with respect to $N$, then its positive equidistants $V_s$, $s >0$  are at least as convex as $V$:
$\phi_s ^* \Pi_s \ge \Pi_0 \ge 0$ for $s > 0$.
\end{proposition}

\vskip .3cm

{\sc End of proof for the 2nd term.}  By the O'neill formula for curvature \cite{Oneill} (see Cor. 1, eq (3), p. 466), the base space $B$ of a
Riemannian submersion  is non-negatively curved provided its total space $Q$ has zero (or positive) curvature.
Applying this to $\pi: M(d,d) \to  Sh$ we get that $Sh$ is   a non-negatively curved manifold at all smooth points. 
(Indeed the sectional curvature of a  two-plane in $T_p Sh$ which is spanned by \on vectors $v, w \in T_p Sh$ is $\sigma = \frac{3}{4} \| F_p (v,w) \|^2$ 
where $F$ is curvature of the \Ri submersion when viewed as a principal $SO(d)$-bundle.  Lemma 2, p.461 of \cite{Oneill}
describes  the relation between the  $A$ occuring there within cor. 1, eq (3), ant the   curvature.) 

By  proposition \ref{curv}  and the fact that $\Sigma= \Sigma_0$ is totally geodesic
at each   smooth point, 
we have  that  each  $\Sigma_s, s > 0$ is concave relative to $\nabla S$, which is to say
that $Q_p \le 0$ for $S > 0$.  It follows that $H \ge 0$ and by symmetry,
as above, we have that $Q = - S H$ with $H$ a positive semi-definite form,
as desired. 

QED for the second term and the proof of proposition \ref{basic}.

\begin{figure}[ht]
\scalebox{0.2}{\includegraphics{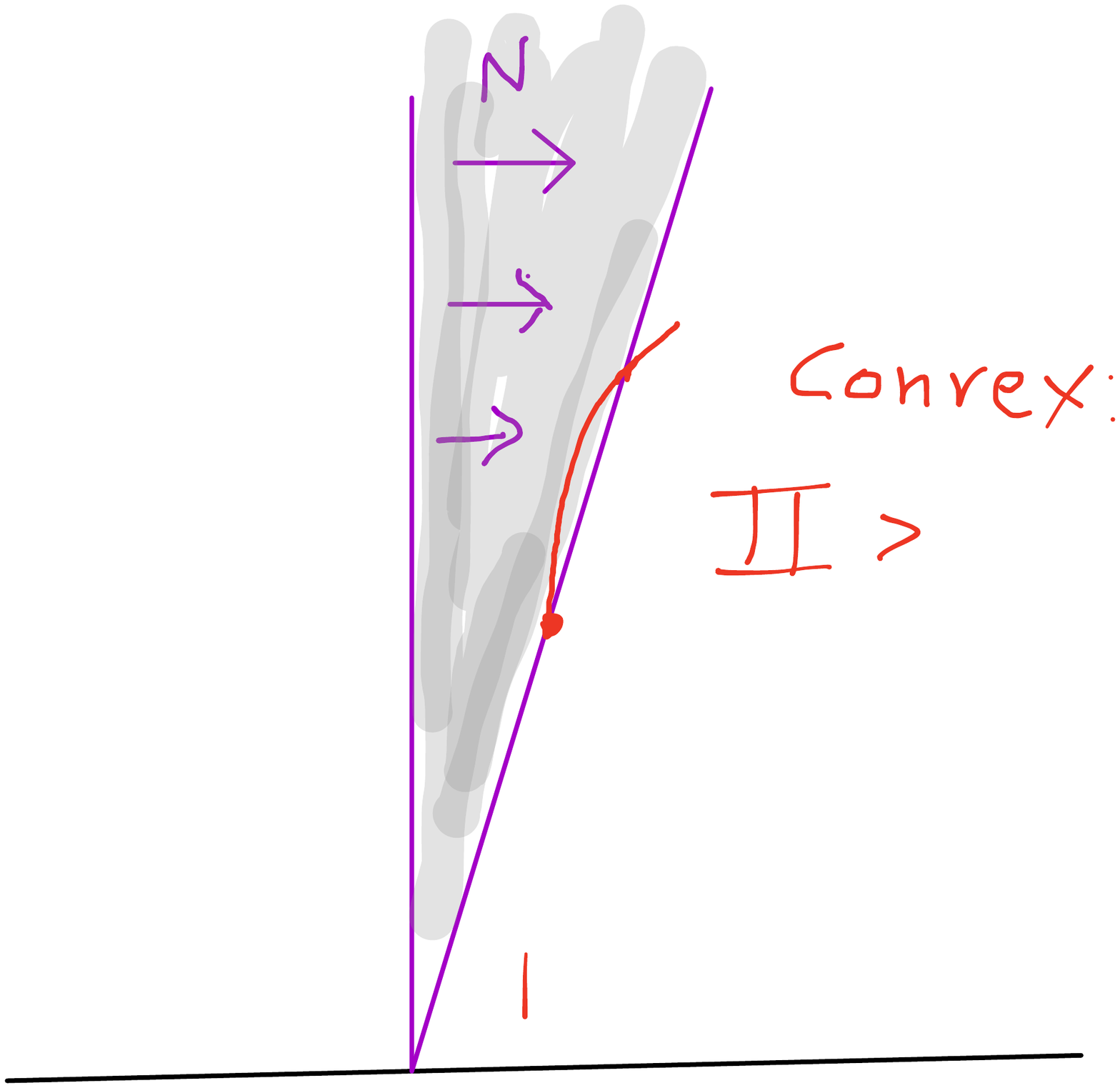}}
\scalebox{0.2}{\includegraphics{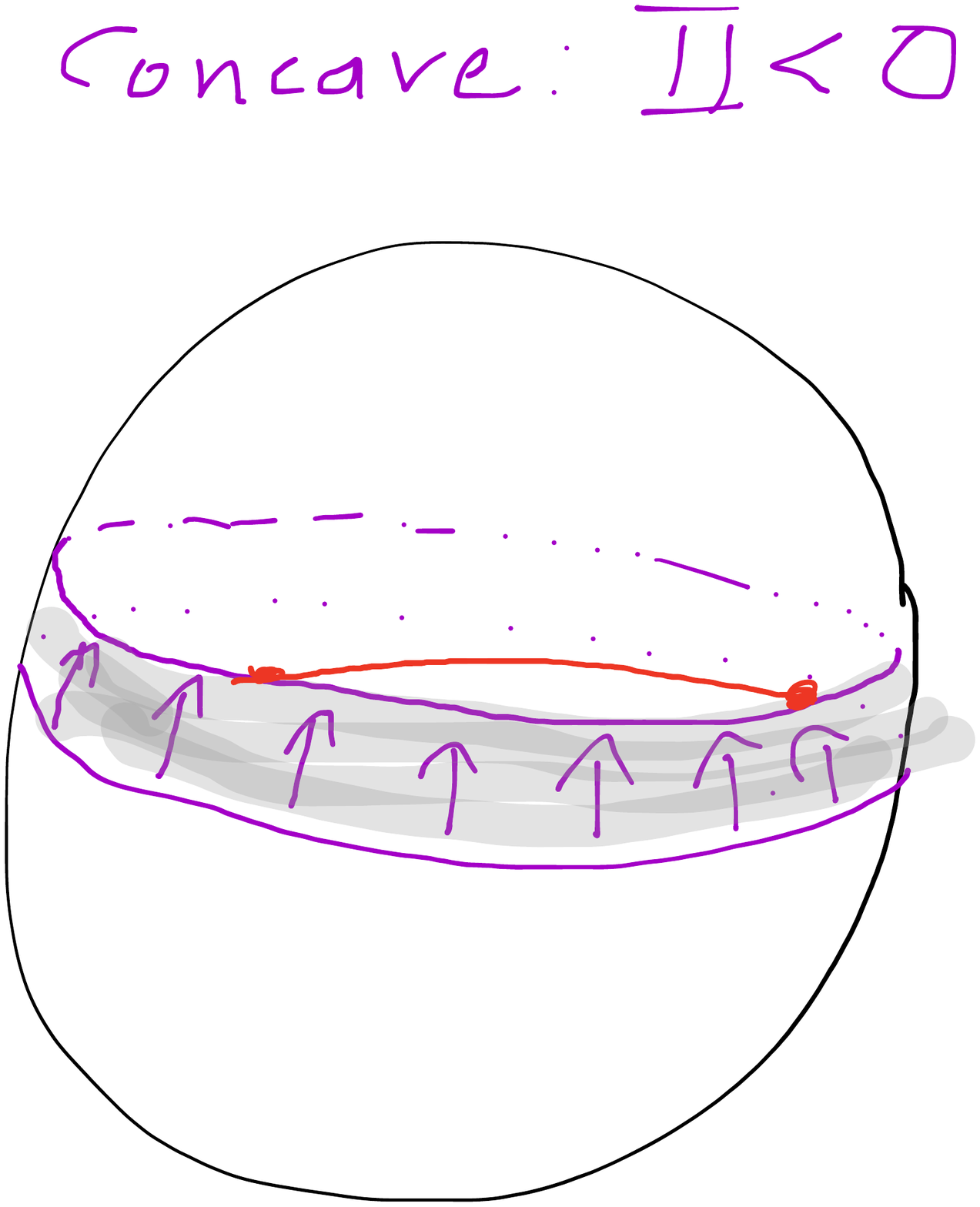}}
\caption{The relation between the sign of the sectional curvatures
and convexity of equidistant hypersurfaces to a totally geodesic submanifold.
The left figure depicts an equidistant from a geodesic in the hyperbolic plane (ambient   curvature  $-1$.
The right figure pictures an equidistant from a geodesic on   the sphere (ambient intrinsic  curvature  $1$).  The first
equidistant is convex relative to the normal while the second is concave.}
\label{fig:indicatrix} 
\end{figure}
 
 \section{Proofs of Theorems}
 
We prove   theorems  \ref{thmA},  \ref{thmB} and \ref{thmC} by strengthening proposition \ref{basic}:  
\begin{proposition} \label{S_convexity} Regardless of whether or not $S$ is smooth along the zero angular solution  $\gamma$ to Newton's equations, 
the composition $S \circ \gamma$  is a convex function of $t$ for  $S >  0$ and  a concave function  for $S <  0$.
If $\gamma$ is bounded
with bounds $r_{ab} \le c$ then $S \circ \gamma (t) = 0$  for at least one time $t$   in each   time interval
of length $\Delta t = \pi (c^3/GM)^{1/2}$, in the Newtonian potential case  and length  $\pi (GM \delta)^{-1/2}$ for the
general potential case as per hypothesis  (\ref{hyp1}), (\ref{hyp2}) and (\ref{hyp3}). 
\end{proposition}
 
{\sc Proof of Proposition \ref{S_convexity} .} 

We first consider the case when $S$ is smooth along $\gamma$, treating the general case as a limit of the smooth case.

If $S$ is smooth along  $\gamma$ then  proposition \ref{basic} asserts that   $\ddot S = -S g$ with $g >  0$ and smooth.
The convex/concave properties of $S \circ \gamma$ follow immediately. 
In case the bounds on the $r_{ab}$ are in force then we know that   
and $g \ge \omega^2 = GM/\delta$ with $\delta$ as per
hypothesis (\ref{hyp1}), (\ref{hyp2}) and  (\ref{hyp3}) in the case of general two-body potential and $\delta = 1/c^3$ in the particular case of the  Newtonian potential.
Compare our differential equation for $S$ to the oscillator equation $\ddot S = - S \omega^2$. The solutions
of the later,  being $S = A sin(\omega (t- t_0))$, have successive   zeros $t_0, t_1, \ldots$ spaced regularly 
at increments  of length  $\pi/\omega$.  By
the   Sturm comparison theorem,  between any two of these zeros lies a zero of our $S$.  Since $1/\omega = \sqrt{ \delta/GM}$
this yields the result for the smooth case. 
 
For the general case, it will suffice  to know that set of   points at which    $S$ fails to be smooth has codimension $2$.
We call this the singular set of $S$ and denote it by $Sing(S)$.   
This assertion regarding the codimension of  $Sing(S)$ is   theorem \ref{SSSVD} of  the next section.  
 
Assuming the validity of this codimension  theorem  \ref{SSSVD},   let $\gamma$ be a zero angular momentum solution to Newton's equations.  Then, by using the
smooth dependence of solutions on initial conditions,  we can find
a family of solutions $\gamma_{\epsilon}$ in $M(d,d)$
which avoids  $Sing(S)$ and converges in the uniform ($C^0$)  topology (or even  $C^k$ topology for any $k$)  to  $\gamma$
on compact time intervals as $\epsilon \to 0$.    By lemma 1,   each   $S \circ \gamma_{\epsilon}$ is convex wherever  it is positive
and concave wherever negative.  The properties of being convex or concave are  closed in the $C^0$-topology, i.e.
the uniform limit of convex functions is convex.   Since $S \circ \gamma_{\epsilon} \to S \circ \gamma$ in the $C^0$ topology our result
for  the convexity / concavity of $S \circ \gamma$  follows.   We proceed to the boundedness implications.
  If the original $\gamma$ satisfies  $r_{ab} (\gamma(t)) \le c$, then its approximating curves   $\gamma_{\epsilon}$ almost satisfy this bound,
  namely, they satisfy $r_{ab} (\gamma_{\epsilon} (t)) \le c + o (1)$ as $\epsilon \to 0$,
  since they $C^0$-converge to $\gamma$.  Thus, by the preceding paragraph,  each  
 $S \circ \gamma_{\epsilon}$ has a  zero in any interval of length $\Delta t = \pi /\omega( \epsilon)$ with $\omega(\epsilon)
 =  GM/\delta(\epsilon)$ and  $\delta(\epsilon)  = \delta (c + o(1))$ as per eq (\ref{hyp3}) above. (Explicitly for  the Newtonian case, 
 $\delta(c + o(1)) = 1/ (c + o(1))^3$.)  Since  $S \circ \gamma_{\epsilon} \to S \circ \gamma$
 we see that    $S \circ \gamma$ must have  a zero in every time interval of length $\pi/ \omega = \pi \sqrt{ \delta/GM}$.   
 
 QED
 
 {\sc Proof of Theorems.} 
 Theorems \ref{thmB} and \ref{thmC} follow immediately from proposition \ref{S_convexity}.
  Theorem \ref{thmA} is the case $d=3$ of theorem \ref{thmB}.
  
  \vskip .3cm
  All that  remains to do now in the way of proofs is to  establish that the codimension of $Sing(S)$ is 2.  

\section{Singular set of the signed distance and the Singular Value Decomposition. }
 
In this section we compute the codimension of   $Sing(S) \subset M(d,d)$ 
 (theorem \ref{SSSVD}).  
 
 
 \subsection{Democracy group. SVD}

 The  signed distance function $S$ enjoys a larger symmetry group than  Newton's equations.
 These additional   symmetries form    the ``democracy group'' and are crucial to  
 identifying  $Sing(S)$.  
   
We saw in subsection \ref{generalD} that the translation-reduced  configuration space is the space of square
matrices  $M(d,d)$, that $\Sigma \subset M(d,d)$ is given by $det(q) =0$
and that by choosing an appropriate basis for ``mass label space'' we can insure that the mass inner product
agrees with the standard Euclidean inner product sot that the    norm squared  of a matrix is  $tr (q^t q) = \Sigma_{i,j} q_{ij}^2$.  
 By inspection, the action of $O(d) \times O(d)$ on $M(d,d)$ by 
 \begin{equation}
 \label{action}
 q \mapsto g_1 q g_2 ^t, g_1, g_2 \in O(d).
 \end{equation} 
is an isometric action which   preserves $\Sigma$. 
It follows that $|S(q)|$,  the distance from $q$ to $\Sigma$ is   invariant under this 
group action.    The action does not quite preserve our  signed distance, since the $O(d)$'s can reverse
orientation.   Indeed
$$S(g_1 q g_2) = \pm S(q); \pm = det(g_1) det(g_1).$$
from which it follows that the action of 
 $SO(d) \times SO(d)$ preserves $S$.  

{\sc Terminology: Democracy group.}   The first, or left $O(d)$ action ($g_1$,  in eq (\ref{action}))  is the usual action of rotations and reflections.
  The second, or   right  $O(d)$ ($g_2$,  in eq (\ref{action}))  is called  the ``democracy group'' since its
  action on the matrix space corresponds to choosing new basis for the mass label space,
  so  in   essence,  permutes, or
 ``democratizes''' the mass labels.  
 
 The  Singular  Value Decomposition [SVD] from Matrix theory \cite{Stewart} 
is a normal form theorem for this group action (\ref{action}). 
 This  decomposition asserts that for any $q \in M(d,d)$
 there is a  {\it diagonal matrix } $x$ and  matrices $g_1, g_2 \in O(d)$
 such that 
 \begin{equation}
 \label{SVD1}q = g_1 x g_2 ^t,  \hskip 1cm [SVD1]
 \end{equation}
Moreover the $g_i$ can be chosen so as to force every nonzero entry of $x$ to be positive,
and the   diagonal entries to be listed  in descending order, thus: 
 \begin{equation}
 \label{SVD2}x = diag(x_1, x_2,\ldots , x_d),  x_1 \ge x_2 \ge \ldots \ge x_d \ge 0.  \hskip 1cm [SVD2]
 \end{equation} 
The diagonal  $x$ written in this form is unique.  Its diagonal entries   $x_i$ are called the ``$i$th principal values'' of $q$. 
The $x_i ^2$ are the eigenvalues of  both of  the symmetric operators $q^t q$ and of $q q^t$. 
 \begin{proposition}
 \label{propSVD} The distance function $|S(q)|$ of $q \in M(d,d) $ to $\Sigma$ is
 equal to $x_d$ above, the $d$th (smallest) principal value of $q$.
 \end{proposition}
 We prove this proposition in the next subsection, below. 
 
 If we impose the constraint that $(g_1, g_2) \in SO(d) \times SO(d)$ when performing the
 normal form computations,  then we get the following
 `specialized' version of the SVD  called the ``pseudo-singular value decomposition'' by \cite{Le} (see p. 361).  
 \begin{proposition} [PsSVD]  
 \label{propSSVD}Given any $q \in M(d,d)$ there is a pair $(g_1, g_2) \in SO(d) \times SO(d))$ and a unique
 diagonal $x= diag(x_1, x_2, \ldots, x_d)$ satisfying $x_1 \ge x_2 \ge \ldots \ge x_{d-1} \ge |x_d|$
 such that 
 $$q = g_1 x g_2 , g_i \in SO(d).$$
 Then 
  $$S(q) = x_d$$
and  $sign(x_d) = sign(det(q)) = sign(S(q))$.  
  \end{proposition}
  In words, the signed distance $S$  is the last `signed' singular value of $q$ in the pseudo-singular value decomposition. 
  
  {\sc Proof of Prop \ref{propSSVD} assuming Prop \ref{propSVD}.}   
  The value of  $det(q)$  cannot be changed by acting on it by
  $(g_1, g_2) \in SO(d) \times SO(d))$ and is equal to $x_1 x_2 \ldots x_d$ if $q = g_1 x g_2 ^t$ with $x = diag(x_1, \ldots , x_d)$.
  Now use the SVD for $q$.   If either one of the elements $g_i $ of the SVD for $q$ is  in $O(d)$ but  not in $SO(d)$ then we can premultiply that element 
   by $diag(1,1, \ldots, 1, -1) \in O(d)$ to get a new $g_i \in SO(d)$ at the expense of perhaps 
  changing $x_d$ to $-x_d$.  Keeping track of the  signs of $det(q)$ and of  $S$ yields  that 
    $S(q) = x_d$, the last  `special' (or `signed')  singular value.
  
  QED   
  
  Finally, here is the assertion we need to  complete all our proofs. 
   \begin{theorem} \label{SSSVD} The   signed distance function $S:M(d,d) \to \R$    is smooth at  any point $q$ of $M(d,d)$ whose smallest two  principal values 
 are distinct.   The complementary set, the singular locus of $S$, is the   set of matrices $q$ whose $d$th and $d-1$st singular values are equal: $x_{d-1} = |x_d|$. 
  This locus is a semi-algebraic set of codimension $2$ within $M(d,d)$.  
    \end{theorem}

  \subsection{A slice.}
  
    The fact   underlying the   proofs of the propositions and  theorems just stated (theorem \ref{SSSVD} etc)  
 is that  the linear subspace $D \subset M(d,d)$ of  diagonal matrices  
 is a {\it  global slice} for our  $O(d) \times O(d)$ action (eq (\ref{action})) on $M(d,d)$. 
 Recall that the {\it orbit} of $q \in M(d,d)$ under this  action is the set 
 $\{ g_1 q g_2 ^t: g_1, g_2 \in O(d) \} \subset M(d,d)$ and that, from basic manifold
 theory, the orbit  is a smooth submanifold. 
 The  assertion that $D$ is a slice for the action  means a number of things
 \begin{itemize}
\item{(a)} every $O(d) \times O(d)$ orbit  intersects $D$
 \item{(b)} the orbit intersects $D$ orthogonally   
 \item{(c)} the intersection is transverse for generic orbit (i.e generic $q$)
 \end{itemize}
 
 Assertion (a) follows from the SVD.
 
Assertion (b)  is a computation.  Let $\xi_1$ and $\xi_2$ be skew symmetric matrices representing
elements of the Lie algebra of our $O(d)$'s, understood to represent the derivatives of the $g_i$ along 
curves passing through $g_i = Id$.    Then the 
tangent space to the orbit through $x$ for   $x \in D$ of the orbit consists of all $d \times d$  matrices $v$  of the form
\begin{equation}
\label{tgtOrbit}
v= \xi_1 x - x \xi_2.
\end{equation}
One sees by direct computation that the diagonal entries of $v$ are all zero, so that $v \perp D$.

Assertion (c) follows  by 
taking ``generic'' matrix to mean one  all of whose principal values are {\bf distinct},
and then making a  more detailed computation  based on the orbit tangent space equation (\ref{tgtOrbit}).  
If we take $\xi_2 = - \xi_1$
in that equation and set $\xi = \xi_1 $ then we compute that $v$ is skew-symmetric with entries    $ (x_i + x_j) \xi_{ij} $ where $\xi_{ij}$ are the entries of $\xi$.
On the other hand, if we take  $\xi_2 = \xi_1 = \xi$ in eq(\ref{tgtOrbit}) we obtain  that  $v$ is a  symmetric matrix 
with  entries   $ (x_j - x_i) \xi_{ij}$.
Now if the $x_i$ are the distinct principal values, we  have that $x_i \pm x_j \ne 0$ for all $i\ne j$  and
it follows easily from this we can obtain any skew-symmetric matrix as a $v$ as per eq (\ref{tgtOrbit}), and  that we can also obtain any symmetric matrix  $v$ which has 
zeros on its diagonal.  Since any matrix at all is   the
sum of a symmetric and a skew-symmetric matrix   we see that the tangent space to the orbit at a generic $x$
consists of all matrices  $v$ with zero entries on the
diagonal, which comprises the orthogonal complement to $D$.

\vskip .3cm

{\sc Proof of Proposition \ref{propSVD}.} As noted just after we introduced the action in eq (\ref{action}), 
the  distance function $|S(q)|$   is invariant under the $O(d) \times O(d)$ action:  
$$|S (g_1 x g_2 ^t)| = |S(x)|.$$
Now  $det(x) = x_1 x_2 \ldots x_d $ so   that 
$\Sigma \cap D = \{ x_1 x_2 \ldots x_d  = 0  \} $   is the union of the $d$ coordinate hyperplanes 
$x_i =0$.  The metric on $M(d,d)$ is Euclidean in the entries,
and $D$ is a $d$-dimensional linear space and in particular totally geodesic: any minimizing
geodesic connecting points of $D$ is a line segment within $D$.  
This implies that for $x \in D$ the $M(d,d)$-distance of $x$ to $\Sigma$  equals the $D$-distance
of $x$ to $\Sigma \cap D$, that is the distance as realized by line segments  {\it within} $D$.
It follows that the problem of computing that distance is a problem in Euclidean geometry.

To solve the problem, let us first fix attention to the case $d=3$.    Observe that  $x_1, x_2, x_3$ are
\on linear coordinates on $D$.   
 The Euclidean  distance of $(x_1,x_2, x_3)$ from the plane $x_1 = 0$ is $|x_1|$.
Since $\Sigma \cap D$ is the union of the three planes $x_1 = 0$, $x_2 = 0$
and $x_3 = 0$, we have that
$$|S(x_1, x_2, x_3)| = min_i  |x_i|.$$
But this minimum {\it} is the 3rd singular value of $x$, namely $x_3$ when
the diagonal values are listed as per the SVD.    The same logic works
for general $d$ and yields $S(x_1, \ldots , x_d) = min_i |x_i|$, which is by
definition the $d$th singular value of $q$.   This proves proposition 
\ref{propSVD}.   

\vskip .3cm 
 
\begin{figure}[ht]
\scalebox{0.2}{\includegraphics{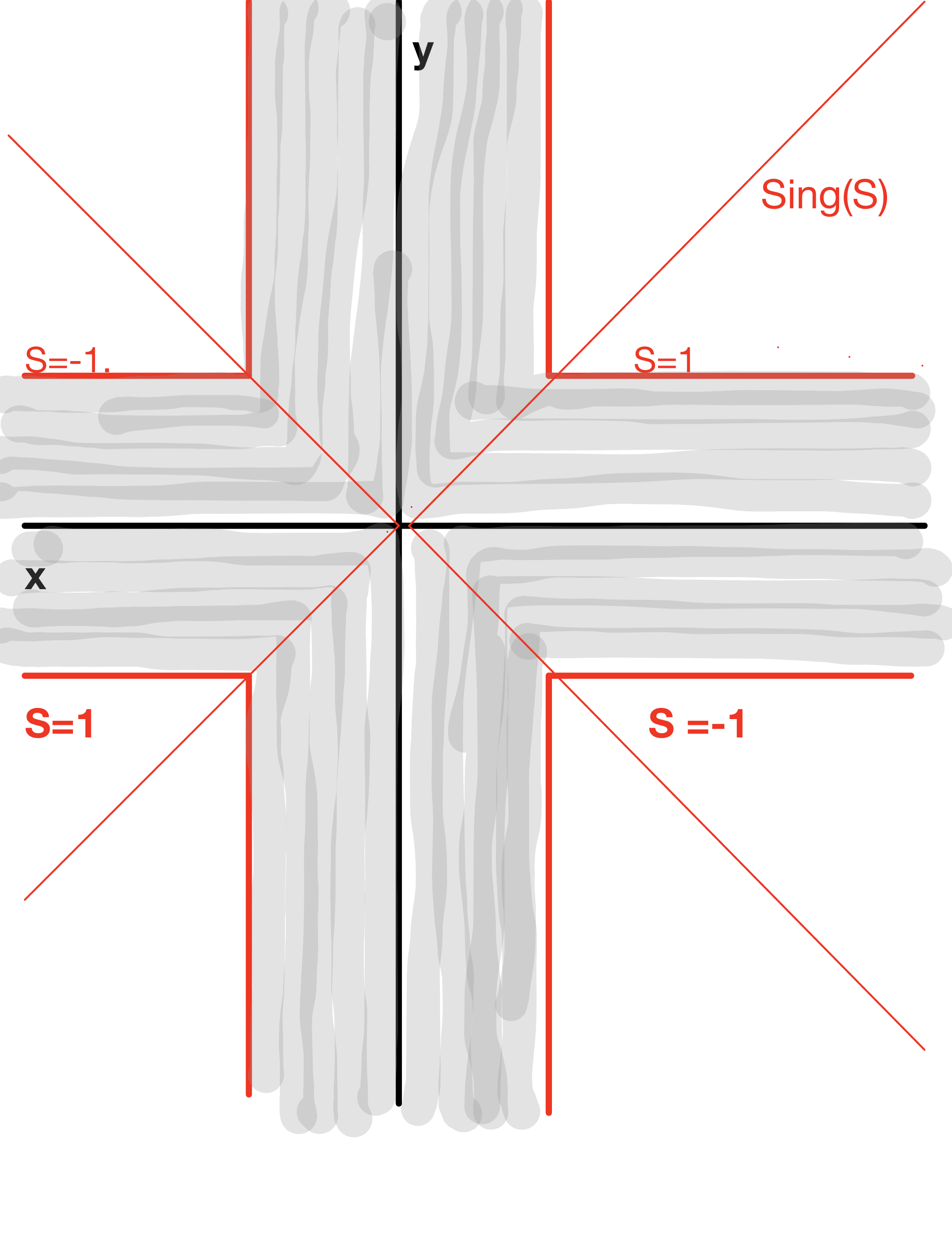}}
\label{cross} 
 \caption{Equidistant curves  to  a cross $xy = 0$
 have corners at which the distance function $|S|$ fails to be smooth. 
 This picture models the contours of $S$ restricted to  the diagonal slice $D$ for $d=2$. The thin red diagonal lines
 indicate $Sing(S) \cap D$.}
\label{fig:indicatrix} 
\end{figure}

{\sc Proof of Theorem \ref{SSSVD}.} 
 {\it The case $d=2$.}  We begin with the case $d=2$ for simplicity and intuition. 
  The configuration space is $M(2,2)$.  The degeneration locus $\Sigma = \{det(q) = 0\}$   is a quadratic cone of signature $(2,2)$ in the vector space  $M(2,2)$.   
 The group $O(2) \times O(2)$ acts isometrically on the matrix space and the diagonal matrices $D$ form a global slice
 as described above. Write    $q = diag(x,y) \in D$.  Then $D \cap \Sigma$
 forms the  ``cross'' $xy=0$.  Within the plane the distance function is 
$S(x,y) = sign(xy) min(|x|,|y|)$.  {\it See figure \ref{cross}. }  The  non-smooth locus of $S$ is the line $x = y$ and $x = -y$
corresponding to the matrices $xI$ and $xJ$ where $I$ is the identity and $J = diag(1,-1)$
It now follows from symmetry that $Sing(S)$ is the union of two two-dimensional conical varieties intersecting
at the origin, namely $\R SO(2) I$ and $\R SO(2) J$.  Taken together this set  is simply $\R O(2)$,
since $J \in O(2)$ and $det(J) = -1$.   If the masses are all equal then this
singular locus corresponds to the Lagrange points (equilateral triangles)
with one cone corresponding to the positively oriented Lagrange configurations
(the north pole of the shape sphere) and the other cone to the negatively oriented
Lagrange configurations.

\eject

{\it The case $d=3$. }  The diagonals are still a slice for the $SO(3) \times SO(3)$ action, and $S$ is invariant under this action.
It follows that we can understand the singularity set of $S$ by looking at its behaviour on the diagonal matrices $diag(x_1, x_2, x_3)$. 
First, suppose we are at a point where  all $x_i > 0$, and such that   $x_1 > x_2 > x_3 > 0$.  Then, $S = x_3$ in a neighborhood of 
our point, which is clearly smooth.
As we move from this point  towards $\Sigma$ along a geodesic orthogonal to $\Sigma$, the value of  $x_3 = S$ steadily decreases
until we hit $x_3 = 0$ at which point $S$  continues to decrease, but smoothly.    
The equality   $S = x_3$ continues into the
region $x_3 < 0$
as long as $x_1 > x_2 > |x_3|$.    This phenomenon is invariant under permutations of the coordinate indices.
Indeed, restricted to   $D$, we have  that $S(x_1,x_2, x_3) = x_i \text{ where } |S(x_1, x_2, x_3)| = |x_i|:= min_k |x_k|$.
Thus the singular locus of $S$ {\it restricted to } $D$ lies on the locus where   $|x_i| = |x_j|$ for some $i \ne j$.
This  locus  is the
 union of 6 planes in $D$, so has dimension $2$, or codimension $1$,  {\it within}  $D$.  
 (The singular locus of the restriction of $S$ to $D$ is a bit smaller that the union of these planes, since we do not need that all three
 principal values are distinct, but only the bottom, two, i.e  we only need $x_2 \ne x_3$ if $x_1 \ge x_2 \ge x_3 \ge 0$ are the singular values.)
 
 At first glance, one guesses that since the singular set  has codimension $1$
 within $D$, then it has  overall  codimension $1$ within $M(3,3)$.  {\it This logic is  wrong.}  Points on the  singular set of  $S$ are {\it not} generic
 with respect  to the $SO(3)\times SO(3)$ action: their symmetry type jumps. Orbits though points of $Sing(S)$ have   dimension  $5$ or less,  not  $6$ like the dimension of  a generic point. (That the orbit through a   generic  point of $D$ is  6-dimensional   is item (b) of `slice' above.)
Since $S$ is invariant under our group $SO(3) \times SO(3)$,
 so is its singular set, $Sing(S)$.   Thus the singular set is the  union of the orbits through  the singular points of the restriction  of  $S$ 
 to $D$.    $Sing(S) \cap D$ has   dimension $2$. If  the  orbit through any point of $Sing(S) \cap D$
 has   dimension $5$ or less then  the singular set itself
has  dimension at most  $7 = 2 +5$. Our space $M(3,3)$ has   dimension $9$, which yields the claimed codimension of $2$.

 It remains to establish that the  orbits through points $x \in Sing(S) \cap D$ have  dimension $5$ or less.  
 The dimension of an orbit of Lie group action is the dimension of the group minus
 the dimension of the isotropy subgroup of that point.  Our group has  $6$.  We show that
 the isotropy group at such a point $x$ has dimension  at least one.  Write $x = diag(x_1, \lambda, \lambda)$ for such a singular point. Let $g(t)$ be  the rotation about the 1st axis by $t$ radians, and $g(-t)$ its inverse.  Clearly $g(t) x g(-t) = x$, establishing that
 the isotropy group is {\it at least} one-dimensional, and hence the orbit has dimension $5 = 6-1$ {\it or less}.
 (A linear algebra computation, following equation (\ref{tgtOrbit}) , shows that this dimension is exactly $5$ as long as $x_1 \ne \lambda$,
 but is unneccessary here since all we need is that the codimension of $Sing(S)$  is at least $2$.) In case $x = diag(x_1, \lambda, - \lambda)$
 with $S(x_1, \lambda, -\lambda) = \lambda$ so that $|x_1|  \ge  | \lambda| \ge 0$, use left multiplication by the matrix $g_1 = diag(-1,1, -1) \in SO(3)$ to replace this
 $x$ by $x = (-x_1, \lambda, \lambda)$  which lies on the same orbit as the original $x$ but now has the
 form of the computation just made.  Since  the orbit is homogeneous its dimension does not
 depend on where on the orbit we choose to compute dimension, and we arrive again at the fact that
 its dimension is  $5$ or less.   
   
 {\it The case $d>3$. }  The proof   is nearly identical to the case $d=3$.   $Sing(S) \cap D$  has codimension $1$,
 being contained in the  union of the hyperplanes where $x_i = \pm x_j$.    At  a generic point of $D$, which is to
 say, off of these hyperplanes,  the $SO(d) \times SO(d)$ action
 is ``almost free'' : the  orbit's dimenison equals that of  $SO(d) \times SO(d)$, as per item (b) of being a slice above.
At a typical point on one of these hyperplanes 
 the isotropy algebra is again one-dimensional , consisting of rotations of  the double eigenvalue plane. (An  ``atypical'' singular point would be one
 for which the   three smallest  singular values are all equal and here the 
 the isotropy algebra has dimension at least $3$.)  Hence the {\it codimension} of $Sing(S)$ is  $1+1$: 
 $1$ for the codimension within $D$ and $1$ for the extra continuous symmetry dimension (isotropy) associated to  each such double ``eigenvalue''  diagonal matrix.  
(Sign discrepancies such as   $x_i = -x_j \ne 0$ are at first bothersome, but the   trick we used in the previous paragraph  of multiplying by an element
 of $SO(d)$ with $\pm 1$'s to change the entries to $x_i = x_j$ works as before. )

QED

\section{Dynamical Vistas and Open Questions}

{\sc Planar precursor.}   

\begin{figure}[ht]
\scalebox{0.3}{\includegraphics{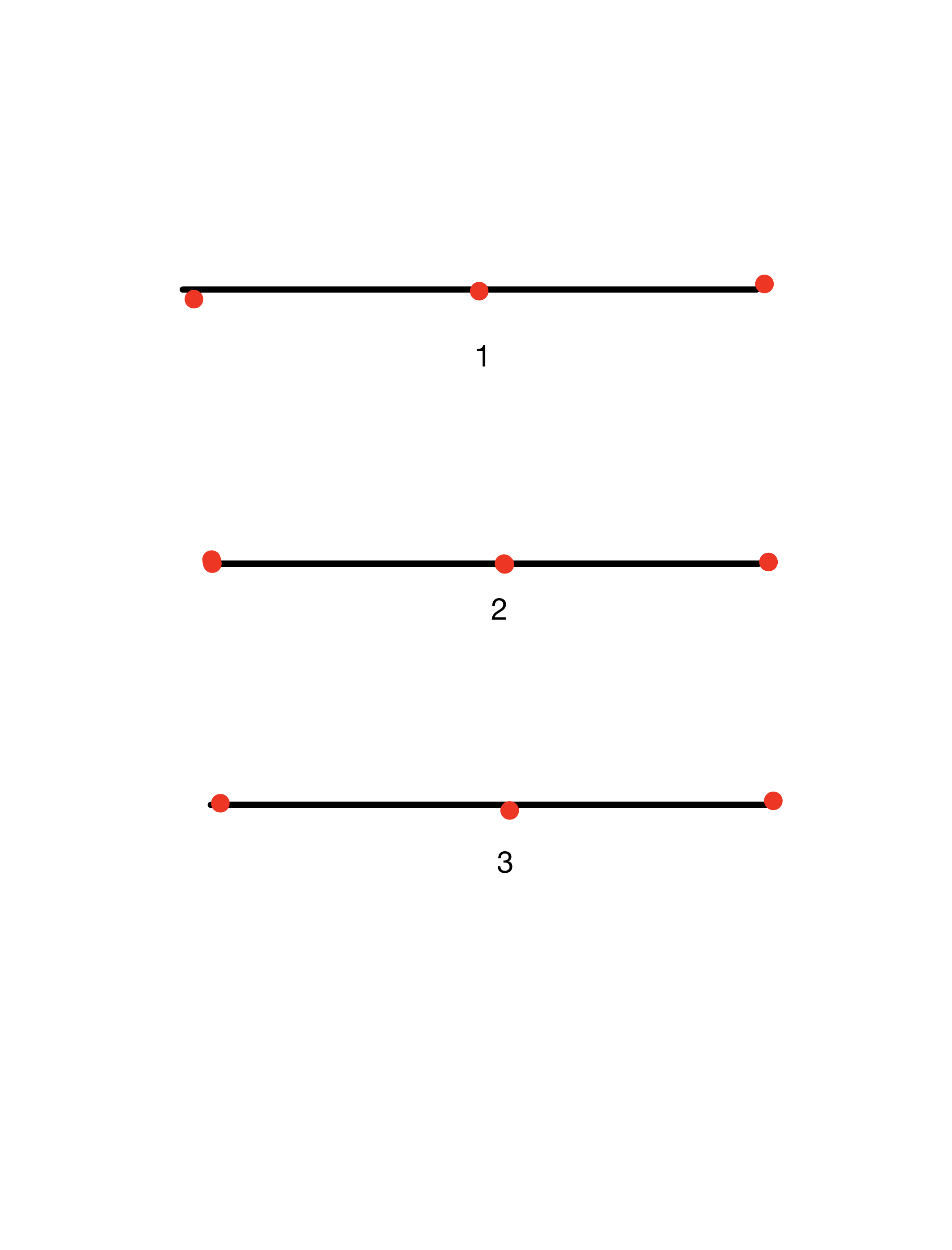}}
\caption{The 3 types of generic collinear 3 body  shapes.}
\label{syz} 
\end{figure}

The planar case of theorem \ref{thmB} or \ref{thmC} asserts that any bounded solution to the planar three-body problem
defined on the whole time line will suffer infinitely colinear instants.  
Colinear instances   are also  called ``syzygies''.   Non-collision syzygies come in three flavors,  1, 2, and 3, depending
on the   mass in the middle. See figure \ref{syz}.  We can thus  associate a syzygy sequence
to such  a solution.  {\it What syzygy sequences are realized? } This question, still largely open, has
motivated much   work.  See for example  \cite{wMoeckel}, and also the closely related
work in which braids (equivalent to ``stutter-reduced'' syzygy sequences)  rather than syzygy sequences are used for the 
symbolic encoding  \cite{Moore}, \cite{Suvakov}, \cite{LiLiao} and references therein.   

\begin{figure}[ht]
\scalebox{0.8}{\includegraphics{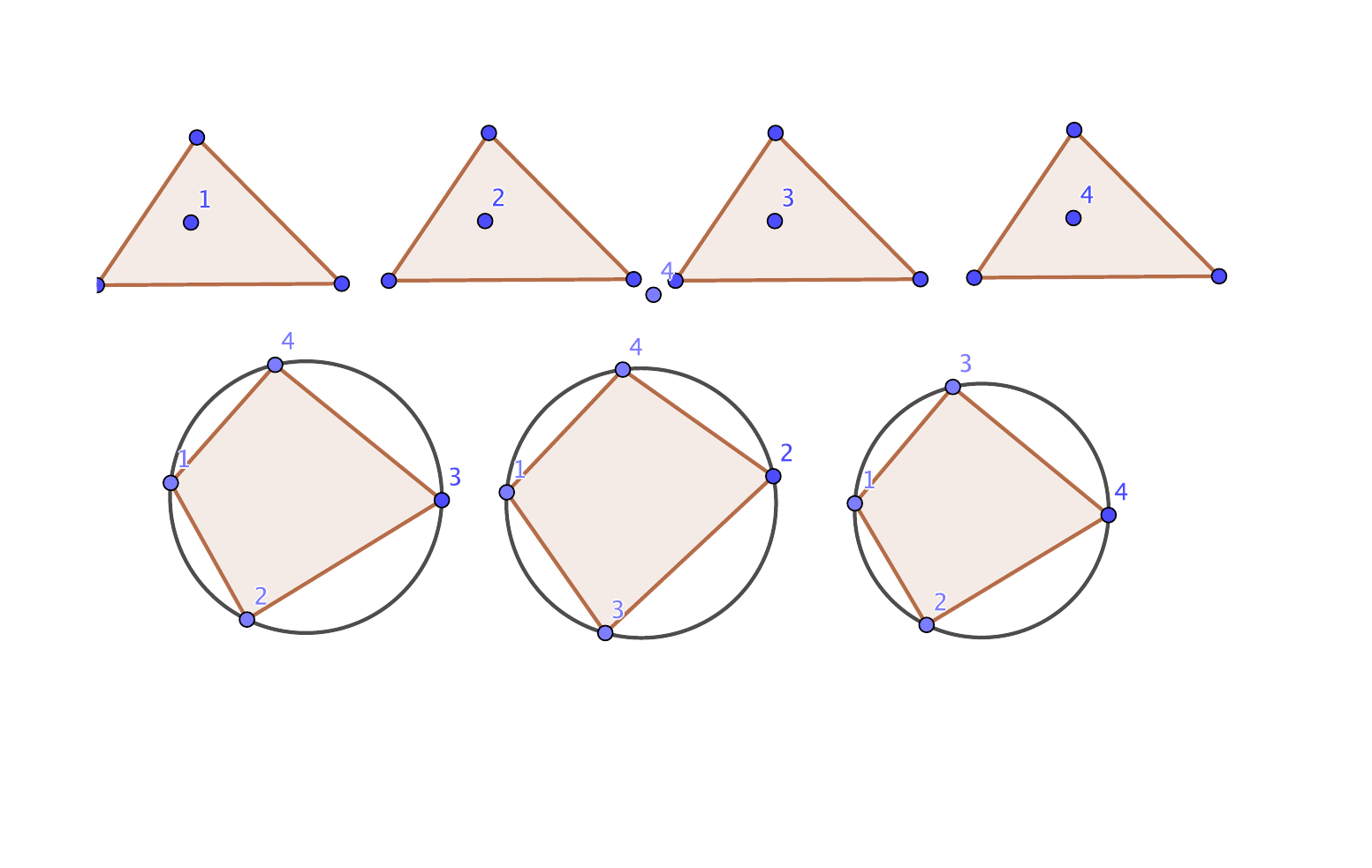}}
\caption{The 7 types of generic planar 4 body  shapes}.
\label{7shapes} 
\end{figure}

Theorem \ref{thmA}   asserts  that  that any bounded solution to the spatial  four-body problem
defined on the whole time line will suffer infinitely coplanar instants.     The generic coplanar configurations divide   into 7   types as per  FIGURE 
\ref{7shapes}.   (We have excluded  as ``non-generic''     configurations for which three of the masses  are collinear.  Binary collision configurations
are thus excluded.)
We now have a   seven letter  alphabet for potential symbol sequences, in analogy with the syzygy sequences of planar three-body dynamics.

{\sc Open questions for the four-body problem in space. }

Q1. Are all possible symbol sequences in this 7-letter alphabet realized by a bounded solution having zero angular momentum?

{\sc Energy and angular momentum considerations}
 Bounded solutions for the Newtonian N-body problem necessarily
 have negative energy. (As soon as  $N >2$ there are negative energy solutions which
 are unbounded.)  Hence the following theorem (see \cite{OnlyLagrange}) represents a
 strengthening of \ref{thmC}for the case $d=2$.   {\bf Theorem:} every zero angular momentum  negative energy solution to the   planar three-body problem 
which does not end and begin in triple collision hits the collinear locus infinitely often. 

We do not  know a single bounded or negative energy solution of the 4-body problem in space
which never suffers co-planarities.  

{\sc More Questions.} 

Q2. Do there exist  negative energy  collision-free solutions of  the spatial four body problem which 
are defined over the whole time line and which  never suffer coplanar instants ?

Q3. If the answer to Q2  is `yes' then are any of these never-coplanar solutions bounded?  

Q4.  If the answer to Q2 is `yes' do  any of these never-coplanar solutions have zero angular momentum?  

In regards to these last two questions, Joseph Gerver has pointed out that there are negative energy collision-free solutions which have no coplanar instants
and   are defined over a time ray  $0 \le t < \infty$. These solutions have nonzero angular momentum. Take the rotating Lagrange  equilateral solution for three of the bodies.
Now take the 4th body to be moving away from this  triple  along the line perpindiular to the plane of the rotating triangle.
If the three masses are equal then the situation is symmetric and the 4th body will stay on this orthogonal  line as it moves out.
The energy of the bound  triple can be taken to be sufficiently small so that the overall energy is negative while
the 4th escaping body escapes hyperbolically to infinity. 
  
\begin{appendix}

 \section{Dispositions: Translation reduction done right.}
  
  We    follow  the idea of ``dispositions''
described  in Albouy-Chenciner \cite{Albouy} in order to identify  the translation  quotient of   configuration space  with the matrix space $M(d,d)$ in a way which is independent of mass choices. We  mostly  stick to the case $d=3$.
 Then we introduce masses and work out how they yield an inner product,  choice of Jacobi vectors
 and  the standard mass-dependent embedding of $M(d,d)$
 back into configuration space.  We strongly recommend \cite{Moeckel}, pp. 34-40 
 for a down-to-earth perspective on this subject. 
 
 As per eq (\ref{config_rep}), we view    the four position vectors $q_a$, $a=1,2,3,4$,  describing
 the instantaneous positions of the four masses as   column vectors and place them side by side to form a    $3 \times 4$
 which we can view as a   linear operator,
$$q  = \left(
\begin{array}{cccc}
 q_1 & q_2  &  q_3 & q_4 
\end{array}
\right) : \R^4   \to \R^3$$
 from  the 
  ``mass label space'' $\R^4$ to our  inertial space of motions $\R^3$.  Thus  the $a$th mass is
  located at the point
$$q_a = q (e_a)$$
where $e_a$, $a =1,2,3,4$  is the standard basis of $\R^4$.
Written out in tensor language,  
$$q = \Sigma q_a \otimes \theta^a$$
where   $\theta^a$ is the basis dual to $e_a$. 

We have identified the  full configuration space with  
 $$M(3,4):=  Hom(\R^4, \R^3) \cong \R^3\otimes \R^{4 *}  .$$  
A translation by $b \in \R^3$ acts on the position vectors $q_a$ as per
the matrix representation of eq (\ref{transl})  which   can be summarized in
tensor language as
$$q \mapsto q + b \otimes \Theta \text{ where } \Theta = \Sigma_a \theta^a$$
since, in matrix terms
 $\Theta = (1,1,1,1)$ so that $b \otimes \Theta = b (1,1,1,1)$ is the matrix
 all of whose columns are $b$.  
 Now $b \otimes \Theta$ is zero on the 3-dimensional linear hyperplane:
  $$L := \ker(\Theta) = \{ \xi \in \R^4: \Theta(\xi) := \xi_1 + \xi_2 + \xi_3 + \xi_4 = 0 \} \subset \R^4$$
from which it follows that the restriction of  $q$ to this subspace is
translation invariant.

 \begin{definition} 
 \label{dispositions} The translation reduction of $q \in M(3,4)$ is its restriction $q_{red}:= q|_{L}: L \to \R^3$,
 to $L \subset \R^4$. 
 The  4-body translation-reduced configuration space is  
 $$Hom(L, \R^3) =  \R^3 \otimes L^*$$
 and the   translation-reduction map is
 $$pr: Hom(\R^4, \R^3) \to Hom(L, \R^3);  pr(q)= q_{red}. $$
 \end{definition}

 {\bf Remark.} Following this prescription we have no   need of masses in order  to identify the quotient of configuration space by translations!   
 
 {\sc \color{red}  Notational Warning. }     Albouy-Chenciner  \cite{Albouy}  use the symbol ${\mathcal D}^*$ for our $L$,
 so that their translation-reduced configuration space is $\R^3 \otimes {\mathcal D} = Hom({\mathcal D}^*, \R^3)$.  Their
 ${\mathcal D}$ is our $\R^{4*}/(\R \Theta) = \R^{4*}/(\R (1,1,1,1)$ which they call the space of dispositions.

  {\sc Degeneration locus.}   The degeneration locus $\Sigma$ coincides with those $q$ for which  $q_{red}:L \to \R^3$ has rank  less than 3.  
 Indeed, a basis for $L$ is formed by $e_{21}, e_{31}, e_{41}$ where $e_{ij} = e_i - e_j$.
 For example, $e_{21} = (-1, 1,0,0)$.  Now $q(e_{ab}) = q_a - q_b$ so relative to this basis, $q_{red}$
 is represented by the 3 by 3 matrix
 $$ \left(
 \begin{array}{cccc}
 q_2 - q_1 & q_3-q_1  &  q_4-q_1  
\end{array}
\right) 
$$
whose determinant is 6 times the signed volume of the tetrahedron whose vertices are  $q_1, q_2, q_3, q_4$.

{\sc Introducing masses.}  Once positive masses $m_a > 0$ are chosen,  we can form the {\it mass vector}
${\vec m} = (m_1, m_2, m_3, m_4) \in \R^4$. Please observe that   
$$\Theta ({\vec m}) = m_1 + m_2 + m_3 + m_4 : = M > 0, $$
demonstrating that  mass vector,  together with the hyperplane $L$,  spans  $\R^4$.  As a consequence 
    $q \in Hom(\R^4, \R^3)$ is uniquely  determined by its  restriction $q_{red}$  to $L$ together with its center
 of mass  :  
 $$q_{cm}:= \frac{1}{\Theta({\vec m})} q ({\vec m} ) = \frac{1}{M} (m_1 q_1 + m_2 q_2 + m_3 q_3 + m_4 q_4) \in \R^3.$$
We   see in this way   that the usual `fixing the center of mass to be zero'' way of forming the translation quotient
 is equivalent to the Albouy  method. 
 

{\sc Mass inner  product, Kinetic Energy, and Jacobi Vectors.} 
The choice of masses also   defines a Euclidean inner product on  $\R^{4*} $  
by declaring that    $\langle \theta^a, \theta^b \rangle = m_a \delta_{a b}$.
In symbols $$ds^2_m = \Sigma m_a (e_a)^2$$
where we are using the fact that the $e_a$ is the  dual basis to   $\theta^a$. 

Let $\V^*, \W$ be real vector spaces.  Choosing  inner products on each  of them
 induces  a canonical  inner product on $ \W \otimes \V^*$ for which 
 $\langle w \otimes \theta, w \otimes \theta \rangle =  \langle w,w \rangle_{\W} \langle \theta,\theta \rangle_{\V^*}  $
In our situation, $\W = \R^3$ comes with its  standard inner product and 
we just used the masses to  put an inner product on  $\R^{4 *}$, so we now  have a mass-dependent  inner product on $M(3,4) = \R^3 \otimes \R^{4*}$.
We call  this inner product the  ``mass inner product'' or sometimes the ``kinetic energy metric''
since it  is the inner product for  which half of  the squared length of a  vector  is the kinetic energy. 
Indeed,  the four instantaneous velocity vectors of   our 4 masses
are obtained  by differentiating $q$ with respect to time
so as to form the  velocity matrix   $\dot q = \Sigma \dot q_a \otimes \theta^a$.
We compute  
$$\frac{1}{2} \langle \dot q , \dot q \rangle  = \frac{1}{2} \Sigma m_a | \dot q_a |^2 = K (\dot q)$$
 which is the   usual expression for the kinetic energy $K$.
 Here we have   reserved the absolute
value sign for the usual norm in our Euclidean $\R^3$,
as in $|q_1|^2 = q_1 \cdot q_1 $. 

The mass inner product  induces an inner product  on our translation-reduced configuration space $Hom(L, \R^3)$ which 
is essential for the body of this paper. There are several equivalent ways to arrive at this inner product.  We will begin  with the
   isomorphism  $\R^{4*} \to \R^{4 **} = \R^{4}$ induced by the mass  inner product on $\R^{4*}$. 
The isomorphism  sends   $\theta^a \mapsto m_a e_a$, $a =1,2,3,4$ and   induces an inner product on $\R^4$ for which
$$\langle e_a, e_b \rangle = \frac{1}{m_a} \delta_{ab}.$$ 
 The isomorphism sends  our
  basic ``translation covector''  $\Theta = \Sigma \theta^a \in \R^{4*}$   to the mass vector ${\vec m}$.
 Since $\Theta$  annihilates $L$ it follows that the mass vector  ${\vec m}$  is   orthogonal to $L$. 
 Thus we get the orthogonal decomposition
$$\R^4 = L \oplus \R {\vec m}.$$
We compute   $\|{\vec m}\|^2 = \|\Theta\| ^2 = \Theta (\vec m) = M$, the total mass.
Now choose an \on basis $E_1, E_2, E_3$ for $L \subset \R^4$, and complete it by adding in  ${\vec m}$ to
form the orthogonal basis $E_1, E_2, E_3 , {\vec m}$ for $\R^4$.
Write the  associated dual  basis for $\R^{4*}$  as  $\omega^1, \omega^2, \omega^3,  \Theta/M$.  Note
that $\omega^i ({\vec m}) =0$. Then the inner product on
$L^*$ can be defined by insisting that the restrictions of the  $\omega^i$   to $L$
forms an \on basis for $L^*$.   Since $Hom(L, \R^3) = \R^3 \otimes L^*$ this
defines an inner product, as desired.

\begin{definition}
\label{JacobiVectors} By a {\it choice of  Jacobi vectors} we mean either  an orthogonal  basis $V_1, V_2, V_3$
for $L$ relative to the mass inner product, or the image vectors $q(V_i)$ of this basis under $q \in Hom(L, \R^3)$.
By a   choice of normalized Jacobi vectors we mean either  an \on basis $E_1, E_2, E_3$  for $L$ or the image of this basis under $q \in Hom(L, \R^3)$. 
These  Jacobi vectors are said to  form an  oriented  basis if the orientation they induce agrees with
that induced on $L$ by the standard basis for $\R^4$ together with the mass vector ${\vec m}$.
 \end{definition}

To better understand the induced  inner product on $Hom(L, \R^3)$, observe
that given  any basis   $u_a$ whatsoever for $\R^4$, and its corresponding dual basis $\omega^a$ for $\R^{4*}$
we can expand out any $q \in M(3,4) $ as $q = \Sigma q(u_a) \otimes \omega^a$. 
If the basis is an orthogonal one, i.e. a triple of `Jacobi vectors' as per the definition
above, then the terms of this expansion are orthogonal relative to our mass metric, meaning that ,
$\|q \|^2 = \Sigma |q (u_a)|^2 \langle \omega^a, \omega^a \rangle$.  
Applying these considerations to our \on basis $E_1, E_2, E_3$ (or `normalized Jacobi vectors'')   we find that 
$$ \| q \|^2 = |X_1|^2 + |X_2 |^2 + |X_3|^2 + M |q_{cm}|^2, \text{ with }  X_i = q (E_i), E_i \text{  \on   for } L.$$
\begin{exer} 
\label{reduced_metric} With $q$ and $E_i$ as above show that the norm squared of the translation reduction $q_{red}$ of $q$
satisfies  
 $$ \| q_{red}  \|^2 = |X_1|^2 + |X_2 |^2 + |X_3|^2 , \text{ with }  X_i = q_{red}  (E_i).$$
 \end{exer}  
 We go a bit deeper into the  last identity for the inner product after some more examples
 around Jacobi vectors

\begin{example} [3 bodies in the plane.] To understand this definition we retreat to the case of 3 bodies in $\R^2$ where Jacobi vectors are better known.
The configuration space is now $Hom(\R^3, \R^2)$ where $\R^3$ is the mass label space and $\R^2$ represents the
inertial space in which the bodies move. 
Then $e_{12} = (1,-1, 0) \in L$ and  $q(e_{12}) = q_1 - q_2$,  is the standard choice of first Jacobi vector.
Take the mass vector to be ${\vec m} = (m_1, m_2, m_3) \in \R^3$ so that the associated
mass inner product on the   mass label space $\R^3$ is given by   $\langle e_a, e_a \rangle = 1/m_a$, $a =1,2,3$.
Then $\langle e_{12}, e_{12} \rangle = \frac{1}{m_1} + \frac{1}{m_2} : = \mu_1 ^2$.
If we write   $m_{12} = m_1 + m_2$ then  $V_2 = (-m_1, -m_2, m_{12}) \in L$. Compute that  
$\langle V_2 , e_{12} \rangle = 0$.   Scale $V_2$ to 
 $U_2 = \frac{1}{m_{12}}(V_1) = (-m_1/ m_{12}, -m_2/m_{12}, 1)$.
Then $q(U_2) = q_2 - (m_1 q_1 +m_2 q_2)/m_{12}$   is the standard expression for the  second Jacobi vector,
and $\langle U_2, U_2 \rangle = \frac{1}{m_{12}} + \frac{1}{m_3} = \mu_2 ^2$
is its normalization factor.  
The vectors 
$E_1 =  \frac{1}{\mu_1} e_{12}$ and $E_2 =  \frac{1}{\mu_2} U_2$
are then an \on basis for $L$, so ``normalized Jacobi vectors''. 
In the standard terminology, the  normalized Jacobi vectors are  $q(E_1) = Z_1$ and $q(E_2) = Z_2$, 
and they yield the identity 
 $|q_{red}|^2 = |Z_1|^2 + |Z_2|^2$.  After identifying $\R^2$ with $\C$, 
 we have $(Z_1, Z_2) \in \C^2$ and this coordinatization of the reduced 
 configuration space as $\C^2$ is the first step in forming shape space and the shape sphere.
 See \cite{meTriangles}  
 \end{example}

 \begin{example}.  [4  bodies in 3-space.]  The three vectors 
 $$J_1 = (1,-1,0,0), J_2 =   (0,0, 1,-1), \text{ and } J_3 = (-p_1, -p_2, p_3, p_4)$$
   form a basis for $L \subset \R^4$
 for any choice of positive $p_i$ with $p_1 + p_2 = p_3 + p_4 = 1$.    Now suppose masses $m_1, m_2, m_3, m_4$
 are given and set $m_{12} = m_1 + m_2, m_{34} = m_3  + m_4,  p_1 = m_1/ m_{12},  p_2 = m_2/m_{12}, p_3 = m_3/m_{34}, p_4 = m_4/m_{34}$
 Then $J_1, J_2, J_3$ are Jacobi vectors -- they are mass-orthogonal. Normalizing them we get
 the corresponding normalized Jacobi vectors $E_i = \sqrt{\mu_i} J_i$ with normalization
 factors $\mu_i$ defined by $\frac{1}{\mu_1} = \frac{1}{m_1} +  \frac{1}{m_2},  \frac{1}{\mu_2} = \frac{1}{m_3} +  \frac{1}{m_4}$, and
  $\frac{1}{\mu_3} = \frac{1}{m_{12}} +  \frac{1}{m_{34}}$.  Then, if we set 
  $\rho_i = q(E_i) \in \R^3$, so that, for example, $\rho_1 = \sqrt{\mu_1} (q_1 - q_2)$,
  and if $q_{cm} = 0$, then $\|q \|^2 = | \rho_1 |^2 + | \rho_2 |^2 +| \rho_3 |^2 $.
  Compare (\cite{Littlejohn}, eqs (2.5), (2.6), pp 2036). 
  \end{example}

 {\sc Returning to the Inner product...}
 
  The restriction-to-L map 
 $pr: Hom(\R^4, \R^3) \to Hom(L, \R^3)$  (see Definition \ref{dispositions} ) implements the quotient of
$ Hom(\R^4, \R^3)$ by the `translation subgroup'' $\T = \R^3 \otimes \Theta \subset \R^3 \otimes \R^{4*} = Hom(\R^4, \R^3)$. 
The choice of masses induces an inner product on  $Hom(\R^4, \R^3)$ and thence on $Hom(L, \R^3)$.  
The masses also induce  an   inclusion $Hom(L, \R^3) \to Hom(\R^4, \R^3)$ which is not present without the additional structure of the masses.
This inclusion and the inner product on  $Hom(L, \R^3)$ are so tightly linked that they are  almost the same thing. 

We formulate all this in more intrinsic linear algebra terms.  Suppose we have an inner product space $\E$
endowed with a subspaces $\T$.  Then $\E/\T$ is naturally endowed with an inner product which makes
$\pi: \E \to \E/\T$ a `submetry' as definied above.  Now consider the  orthogonal
complement $\T^{\perp} \subset \E$.  Then the restriction of $\pi$ to $\T^{\perp}$ is an isometry
between $\T^{\perp}$ and $\E/\T$.  Now apply these considerations to the case that 
$\E = Hom(\V, \W)$,  $\E/\T = Hom(\S, \W)$ where $\S \subset \V$ is a linear subspace of $\V$
and where the projection map $\pi$ is the map which sends $q: \V \to \W$ to its restriction $q|_{\S}: \S \to \W$.  
Then $\T \subset Hom(\V, \W)$ consists of those linear operators which are zero on $\S$.
An inner product on $\V$ and $\W$ induces ones on all these spaces.  
\begin{exer} Show that, continuing with the above terminology, $\T^{\perp}$
consists of those linear operations $\V \to \W$ which are identically zero on the orthogonal complement
$\S^{\perp} \subset \V$ to $\S$ and that the induced isometric inclusion $Hom(\S, \W) \to Hom(\V, \W)$
is the map which takes  a linear operator $q_{red} : \S \to \W$ and extends off of $\S$ to obtain a  map $\V \to \W$
 equal to $q_{red}$ on $\S$ and identically  zero on $\S^{\perp} \subset \V$.
\end{exer}

Applying these   considerations  to our situation of $L \subset \R^4$ and $q_{red} \in Hom(L, \R^3)$  we
 see that the mass induced canonical extension $q = i(q_{red})$ of $q_{red}$ is obtained 
by insisting that  $q({\vec m}) = 0$, since ${\vec m}$ is the normal vector to $L$ relative to the mass metric.
But $q({\vec m}) = M q_{cm}$.   {\bf  The   mass-induced inclusion $i: Hom(L, \R^3) \to Hom(\R^4, \R^3)$ 
of the translation-reduced configuration space  into our original configuration space
$Hom(\R^4, \R^3)$ is
obtained by insisting that the center of mass of the extension
$q = i( q_{red})$ is zero.} Thus the process of using inner products to
implement this inclusion  amounts to the usual `center-of-mass zero' prescription for performing
translation reduction.  This  inclusion $i$ an isometry onto its image, and  we arrive 
in this way at a solution to exercise \ref{reduced_metric}  regarding the  
inner product on $Hom(L, \R^3)$.  Comparing the formulae we see that the
prescription for the inner product amounts to the equation that  
$\|q\|^2 = \| q_{red} \|^2$ whenever $q_{cm} = 0$.

 {\sc  Democracy Group; Linear isometries.}  Let us again suppose that $\V$ and $\W$
 are finite-dimensional real inner product spaces.   Then  $\W \otimes \V^* = Hom(\V, \W)$
 inherits a canonical inner product for which the orthogonal  groups
 $O(\V)$ and $O(\W)$ act  isometrically according to
 $$q \mapsto g_1 \circ q \circ g_2 ^{-1} ,  \text{  for } g_1 \in O(\W), g_2 \in O(\V).$$
Specializing to our situation,
 $\V = L \subset \R^4, \W = \R^3$,  so that   the 6-dimensional Lie group
 $O(L) \times O(\R^3) \cong O(3) \times O(3)$ acts isometrically on our translation-reduced
 configuration space $Hom(L, \R^3)$.    The $g_1$ factor acting by left multiplication
 is the usual action on a 4-body configuration by rotation.  The $g_2$ action is a bit more mysterious. 
 It is not a symmetry of Newton's equations.
 
 \begin{definition} The democracy group action is the action of 
 $O(L)$ on $Hom(L, \R^3)$ by $q \mapsto q g_2 ^{-1}$.
 \end{definition}
 
To better understand the democracy group action,
choose an  \on basis for $L$, which is to say normalized Jacobi vectors $E_1, E_2, E_3$.
This choice   induces   an isomorphism $O(L) \cong O(3)$.
 $\R^3$ comes with a standard basis $u_1, u_2, u_3$ so we can just write $O(\R^3) = O(3)$.   Relative to these two bases
 $q_{rel}: L \to \R^3$
 becomes  a  $3 \times 3$ matrix $X$ with  entries  $X_{ij} = u_i \cdot q(E_j)$.  
 In this way we identify
 $$Hom(L, \R^3) = M(3,3) := \text{ the space of all real three-by-three matrices}.$$
 The  mass metric in these coordinates is simply the standard entry-wise
 Euclidean structure:   $\| q_{rel}\|^2 = tr (X^t X) = \Sigma_{i,j} X_{ij}^2$. 
 Our $O(L) \times O(\R^3)$  action is simply standard matrix multiplication by $O(3) \times O(3)$:
 $$X \mapsto g_1 X g_2 ^t.$$
 The degeneracy locus is  
 $$\Sigma = \{X \in M(3,3):  det(X) = 0 \}.$$
 and is mapped to itself by the $O(3) \times O(3)$ action, since $det(g_1 X g_2 ^t) = \pm det(X)$.
  
 The idea expressed by the terminology ``Democracy group''   is that a choice of Jacobi vectors involves selecting out 
 certain  masses to play  special roles. An element of $O(L)$ changes the basis, i.e the Jacobi
 vectors, and hence corresponds to choosing different sets of masses for these roles.
 It permutes the mass labels.    Indeed if all the $m_a$ 
 are equal, then the  permutation group of the mass labels, 
i.e. of the basis vectors $e_a$ for the label space $\R^4$,  forms  a subgroup of $O(L)$.
 We owe the picturesque name `democracy group' to Littlejohn and Reinsch, \cite{Littlejohn}.
 
 
 \vskip .3cm
 
 \section{Unoriented vs Oriented Shape Space}
 
We describe the  oriented and unoriented shape space and the relation between them.
For the general $N$ body problem in $\R^d$ the configuration space is $Hom(\R^N, \R^d)$
and these  shape spaces are the quotient spaces of the configuration space by
the groups $SE(d)$ and $E(d)$ respectively. 
 The  quotient $E(d)/SE(d)$ is the two-element group (with nontrivial element generated by any orientation reversing isometry)
 so that we expect that the natural map from oriented to un-oriented shape space will be a   $2:1$ 
 branched cover, branched over the degeneration locus   $\Sigma$.   
  
 We form the quotient spaces in stages. We saw   in the previous appendix how the  quotient of the configuration space 
  by translations is  $Hom(L, \R^d) \cong M(d,d)$ in case $N= d+1$. It remains to quotient by the linear isometries,
  i.e. the rotations $SO(d)$ and rotations and reflections $O(d)$.  
  Thus oriented shape space becomes  $M(d,d)/SO(d)$ while unoriented shape space becomes
  $M(d,d)/O(d)$, where $g \in O(d)$ or $SO(d)$ acts on $q \in M(d,d)$ by $q \mapsto gq$. 
 
  {\sc Unoriented shape space.}  The map $q \mapsto q^t q \in Sym(d)$ realizes the $O(d)$ quotient.
  Here $Sym(d)$ denotes the space of symmetric $d \times d$ matrices.
  When we say ``realizes the quotient'' we mean that for $q', q \in M(d,d)$
  there exists a $g \in O(d)$ such that $q' = gq$ if and only if $(q')^t q' = q^t q$.
  This fact follows from a basic theorem from representation theory.
  The matrix $q^t q$ is sometimes called the `Gram matrix', being a matrix of inner products of position vectors.
  Any matrix of form $q^t q$ is positive semi-definite, and any positive
  semi-definite matrix $s$ can be expresses as $s = q^t q$ for some $q \in M(d,d)$.
  (Take $q = \sqrt{s} \in Sym(d)$ for example.).  These facts prove
  that  un-oriented shape  space is  the ``positive semi-definite cone'': the  
   closed convex cone of positive semi-definite symmetric matrices within $Sym(d)$.  
   The  boundary of the cone consists of those 
 non-negative symmetric matrices whose rank is not full and thus  corresponds  to the shape projection of  our   degeneration locus $\Sigma$.
 
{\sc The map between.}  The map from oriented shape space to shape space   maps  onto this shape cone,  and forms  a 2:1 cover
 branched along $\Sigma$.    Indeed,   an unoriented nondegenerate simplex  shape 
  has precisely two oriented representative shapes, one having   positive volume, the other  having negative
 volume, while a degenerate shape in $M(d,d)/SO(d)$ has precisely one representative in the unoriented shape space. 
  
 {\sc A bit of topology.} Take two identical copies of a  closed convex cone with nonempty interior in any finite dimensional real vector space.
 Glue one copy to the other along the boundary,  using   the identity map of the boundary as gluing map.    One checks without great difficulty that the
 result is homeomorphic to the original vector space: we have `blown up' or desingularize the boundary.
 These general considerations may serve to convince the reader that oriented shape space is indeed homeomorphic
 to the Euclidean space $Sym(d)$

\end{appendix}


\begin{thebibliography}{99}

\bibitem{Abraham} R. Abraham and J.E. Marsden, {\bf Foundations of Mechanics},  [1978], 2nd ed., Addison-Wesley, Reading, MA.

\bibitem{Albouy} A. Albouy and A. Chenciner, {\it Le probl\'eme des n corps et les distances mutuelles}, Inventiones {\bf 131} (1998), 151-184.

\bibitem{Arnold}  V.I. Arnol'd,  S.M. Gusein-Zade, and A.N. Varchenko, {\bf Singularities of Differentiable Maps, volume I},
 [1985], Birkhauser, (see esp. p.28-30).

\bibitem{Gromov} M. Gromov,  {\bf Sign and   Meaning of Curvature}, Rendiconti del Seminario Matematico e Fisico di Milano
 {\bf 61}, issue 1,    (1991), 9-123., (see esp. p. 34-37).

\bibitem{Kendall} D. G. Kendall, D.  Barden,  T.K. Carne,   and H.  Le, 
 {\bf Shape and shape theory}, 
    {Wiley Series in Probability and Statistics},
  {John Wiley \& Sons, Ltd., Chichester},
      [1999]. 

\bibitem{Le} H. Le, {\em On Geodesics in Euclidean Shape Spaces}, J. London Math. Soc., (2), {\bf 44},  (1991) 360-372.

\bibitem{LiLiao} X. Li and S. Liao, {\em Collisionless periodic orbits in the free-fall three-body problem},  arxiv:1805.07980.

\bibitem{Littlejohn} R. Littlejohn and M. Reinsch, {\em Internal or Shape coordinates in the n-body problem}, Phys.  Rev. A ,  {\bf 52}, no. 3, (1995),2035-2051.

\bibitem{Moeckel} R. Moeckel, {\em Lecture Notes on Central Configurations}, \url{http://www-users.math.umn.edu/~rmoeckel/notes/Notes.html}.

\bibitem{wMoeckel} R. Moeckel and R. Montgomery. {\it Realizing All Reduced Syzygy Sequences in the Planar Three-Body Problem},
  Nonlinearity,  {\bf 28} (2015), 1919-1935. 

\bibitem{isoholonomic} R. Montgomery, {\it The isoholonomic problem and some of its applications}, 
Comm. Math. Phys., {\bf 128}, (1990),   565-592. 

\bibitem{Infinitely} R. Montgomery, {\it Infinitely Many Syzygies Archives for Rational Mechanics}, {\bf 164}, no. 4,  (2002),  311--340.  (DOI) 10.1007/s00205-002-0211-z.

\bibitem{OnlyLagrange} R. Montgomery, {\it The zero angular momentum three-body problem: all but one solution has syzygies}, 
 Erg. Th. and Dyn. Systems, {\bf 27}, no. 6, (2007),   1933-1946.
 
\bibitem{meTriangles} R. Montgomery, {\it The Three-Body Problem and the Shape Sphere},  Amer. Math. Monthly, {\bf 122}, no. 4,(2015), 299-321. 


\bibitem{tour} R. Montgomery, {\bf A tour of sub-Riemannian Geometry},  Mathematical Surveys and Monographs,  {\bf 91},
American Math. Society, Providence, Rhode Island, [2002].   (see esp. p. 204).   

\bibitem{Moore} C. Moore, {\it Braids in Classical Gravity},  Phys. Rev. Lett., {\bf 70}, (1993), 3675-3679.  

\bibitem{Oneill} B. O'Neil, {\it The fundamental equations of a submersion}, Michigan Math. J. , {\bf 13}, no. 4,  (1966)  459--469,
(see esp. Cor. 1, eq (3), p. 466).
doi:10.1307/mmj/1028999604. 



\bibitem{Strang} G. Strang, {\bf Introduction to Linear Algebra},  (2nd.  ed.), [1976],  Wellesley-Cambridge Press,  ISBN 0-9614088-5-5,
(esp. section 3.4, 141-142). 

\bibitem{Suvakov} M. Suvakov and V. Dmitrasinovic,  {\em Three Classes of Newtonian Three-Body Planar Periodic Orbits}, Phys. Rev. Lett., {\bf 110}, (2013),114301-114305.   arxiv:1303.0181.  

\end{thebibliography}
\end{document}